Publicado en:

# Lógica Doble LD y Gráficos Existenciales Gamma-LD
# LD Double Logic and Gamma-LD Existential Graphs

MANUEL SIERRA-ARISTIZÁBAL[1]

**Resumen.** En este trabajo se presenta el sistema deductivo para la *Lógica Proposicional Doble*, LD, y los *Gráficos Existenciales Gamma-LD*. Se prueba, de manera rigurosa y detallada, la consistencia de LD y que los teoremas de LD corresponden exactamente a los gráficos existenciales válidos de Gamma-LD. Cuando se restringe el lenguaje de LD al lenguaje de la *Lógica proposicional Clásica* LC, entonces la restricción asociada a Gamma-LD coincide con los gráficos existenciales válidos del sistema *Alfa* de Charles Sanders Peirce. Resulta que los teoremas de la *Lógica Proposicional Intuicionista* LI son teoremas de LD, además, cuando se restringe el lenguaje de LD al lenguaje de LI, entonces la restricción asociada a Gamma-LD coincide con los gráficos existenciales válidos del sistema *Alfa intuicionista* de Arnold Oostra. Como consecuencia, se infiere que Gamma-LD tiene como casos particulares, los gráficos existenciales Alfa de LC y de LI. Finalmente, en LD se derivan las definiciones de *verdad* y *falsedad* presentadas por Aristóteles, con las cuales, se ilustra la capacidad que tiene LD para solucionar una versión de la *paradoja del mentiroso*, donde LC y LI fracasan.

**Abstract:** In this paper, the deductive system for double propositional logic (DL) and gamma-DL existential graphs are presented. It rigorously proves the consistency of the DL and that the DL theorems correspond exactly to the valid existential graphs of gamma-DL. When the language of DL is restricted to the language of *classical propositional logic* (CL), the restriction associated with gamma-DL coincides with the valid existential graphs of Charles Sanders Peirce's *alpha* system. It turns out that *intuitionistic propositional logic* (IL) theorems are DL theorems; furthermore, when the language of DL is restricted to the language of IL, the restriction associated with gamma-DL coincides with the valid existential graphs of Arnold Oostra's *intuitionistic alpha* system. As a consequence, it is inferred that gamma-LD has as particular cases, the alpha existential graphs of

[1] msierra@eafit.edu.co ORCID: 0000-0001-8157-2577

.

LC and LI. Finally, in DL, the Aristotelian definitions of truth and falsity are derived, with which the ability of DL to solve a versión of the *liar paradox*, where LC and LI fail, is illustrated.

# Presentación

En este trabajo se presenta el sistema deductivo *Lógica Proposicional Doble* LD, como una extensión generalizada de la *Lógica Proposicional Clásica* LC, la cual incluye operadores alternos de afirmación y de negación. A partir de estos operadores, en LD se definen operadores alternos a los conectivos binarios clásicos, resultando que LD también es una extensión generalizada de la *Lógica Proposicional Intuicionista* LI. No sólo se tiene que LC y LI son casos particulares de LD, sino que, además, los conectivos binarios intuicionistas, sorprendentemente, pueden ser definidos en términos de los conectivos clásicos y de los conectivos alternos en LD (característica que se pierde cuando se restringe el lenguaje de LD al de LI.

En este trabajo también se presentan los *gráficos existenciales* Gamma-LD para LD. Se prueba, de manera rigurosa y detallada, la consistencia de Gamma-LD y que los teoremas de LD corresponden exactamente a los gráficos existenciales válidos de Gamma-LD. Cuando se restringe el lenguaje de LD al lenguaje de LC, entonces la restricción asociada a los gráficos Gamma-LD coincide con los gráficos existenciales válidos del sistema Alfa de Charles Sanders Peirce, presentado en [1] y con los gráficos existenciales válidos del sistema Alfa$_0$ de Yuri Poveda, presentado en [2]. Resulta también que los teoremas de LI son teoremas de LD, además, cuando se restringe el lenguaje de LD al lenguaje de LI, entonces la restricción asociada a los gráficos Gamma-LD coincide con los gráficos existenciales válidos del sistema Alfa-Intuicionista de Arnold Oostra, presentado en [3]. Como consecuencia se infiere que, Gamma-LD tiene como casos particulares los gráficos existenciales de LC y de LI.

Al ser interpretado el operador de afirmación alterna como *Verdad Aristotélica*, y el operador de negación alterna como *Falsedad Aristotélica*, en LD se tienen como consecuencia, las definiciones de verdad y falsedad dadas en [4]: "Decir de lo que es que es, y de lo que no es que no es, es lo verdadero; decir de lo que es que no es, y de lo que no es que es, es lo falso". Resultando que estas definiciones son teoremas de LD. Finalmente, con esta lectura de los operadores alternos, se ilustra la capacidad que tienen LD y Gamma-LD para solucionar una versión de la *paradoja del mentiroso*, presentada en [5] y en [6], donde LC y LI fracasan.

# Lenguaje

Se parte de un conjunto infinito enumerable de *Fórmulas Atómicas Clásicas* **FAC**, las cuales se denotan a, b, c, d, … . A partir de FAC se obtiene el conjunto de *Fórmulas Atómicas Alternas* **FAA**, las cuales se denotan $\underline{a}$, $\underline{b}$, $\underline{c}$, $\underline{d}$, … . El conjunto de *Fórmulas Atómicas* **FAT**, se obtiene de la unión de los dos anteriores FAT =FAC∩FAA. Los conectivos de LD se clasifican de la siguiente manera: conectivo unario clásico = {~}, conectivos unarios alternos = {¬, +}, conectivos binarios clásicos = {⊃, •, ∪, ≡} y conectivos binarios alternos = {→, ∧, ∨, ↔}.

El conjunto **FC** de *Fórmulas Clásicas* se define de la siguiente manera:

1. $a \in FAC \Rightarrow a \in FC$.
2. $X \in FC \Rightarrow {\sim}X \in FC$
3. $X,Y \in FC \Rightarrow X \supset Y,\ X \cup Y,\ X \bullet Y,\ X \equiv Y \in FC$.
4. Sólo 1 a 3 determinan FC.

El conjunto **FI** de *Fórmulas Intuicionistas* se define de la siguiente manera:

1. $\underline{a} \in FAA \Rightarrow \underline{a} \in FI$
2. $X \in FI \Rightarrow \neg X \in FI$.
3. $X,Y \in FI \Rightarrow X \rightarrow Y,\ X \vee Y,\ X \wedge Y,\ X \leftrightarrow Y \in FI$.
4. Sólo 1 a 3 determinan FI.

El conjunto **FOR** de *Fórmulas de LD* se define de la siguiente manera:

1. $X \in FC \Rightarrow X \in FOR$.
2. $X \in FI \Rightarrow X \in FOR$.
3. $X,Y \in FOR \Rightarrow X \supset Y,\ X \cup Y,\ X \bullet Y,\ X \equiv Y,\ X \rightarrow Y,\ X \vee Y,\ X \wedge Y,\ X \leftrightarrow Y \in FOR$.
4. $X \in FOR \Rightarrow {\sim}X, \neg X \in FOR$.
5. Sólo 1 a 4 determinan FOR.

El conjunto **FA** de *Fórmulas alternas* se define de la siguiente manera:

$X \in FA \Leftrightarrow X = \neg Y$ con $Y \in FOR$ o $X = \underline{a}$ donde $\underline{a} \in FAA$. Si $X \in FOR$ entonces $\underline{\mathbf{X}}$ indica que $X \in FA$.

Definición 1. (*Operadores de verdad*). En **LD**, a partir de la *Falsedad alterna* $\neg X$, se definen los siguientes conectivos (donde $X,Y \in FOR$):

*Verdad alterna*: $+X = \neg{\sim}X$.

*Disyunción alterna*: $X \vee Y = +(X \cup Y)$.

*Conjunción alterna*: $X \wedge Y = +(X \bullet Y)$.

*Implicación alterna*: $X \rightarrow Y = +(X \supset Y)$.

*Equivalencia alterna*: $X \leftrightarrow Y = +(X \equiv Y)$.

*Refutable*: $-X = {\sim}+X$.

*Satisfacible*: $\otimes X = {\sim}\neg X$.

*Buena fundamentación*: $*X = (\neg X \cup +X)$.

Como consecuencias se tienen las siguientes equivalencias entre los llamados ***operadores de verdad***:

a. $+X \equiv \neg{\sim}X \equiv {\sim}\otimes{\sim}X \equiv {\sim}-X$.

b. $\neg X \equiv +{\sim}X \equiv {\sim}\otimes X \equiv {\sim}-{\sim}X$.

c. $\otimes X \equiv {\sim}+{\sim}X \equiv {\sim}\neg X \equiv -{\sim}X$.

d. $-X \equiv {\sim}+X \equiv {\sim}\neg{\sim}X \equiv \otimes{\sim}X$.

e. $*X \equiv +X \cup \neg X \equiv -X \supset \neg X \equiv \otimes X \supset +X$.

# Sistema deductivo

El *sistema deductivo para* **LD** consta de los siguientes axiomas (donde $X,Y,Z \in FOR$, $\underline{X} \in FA$):

Ax1.1 $X \supset (Y \supset X)$

Ax1.2 $(X \supset (Y \supset Z)) \supset ((X \supset Y) \supset (X \supset Z))$

Ax1.3 $X \supset (X \cup Y)$

Ax1.4 $Y \supset (X \cup Y)$

Ax1.5 $(X \supset Z) \supset ((Y \supset Z) \supset ((X \cup Y) \supset Z))$

Ax1.6 $(X \bullet Y) \supset X$

Ax1.7 (X•Y)⊃Y

Ax1.8 (X⊃Y)⊃((X⊃Z)⊃(X⊃(Y•Z)))

Ax1.9 X⊃(~X⊃Y)

Ax1.10 X∪~X

Ax1.11 (X≡Y)⊃(X⊃Y)

Ax1.12 (X≡Y)⊃(Y⊃X)

Ax1.13 (X⊃Y)⊃[(Y⊃X)⊃(X≡Y)]

Ax2.1 AxMP+: +(X⊃Y)⊃(+X⊃+Y)

Ax2.2 AxR: ¬X⊃~X

Ax2.3 AxT: $\underline{X}$⊃+$\underline{X}$

Ax2.4 Ax+: X∈{Ax1.1, …, Ax2.3} ⇒ +X es un axioma.

Como única *regla de inferencia* se tiene el *Modus Ponens* MP: de X y X⊃Y se infiere Y.

Definición 2. (*Tipos de axiomas*). *Axiomas clásicos* AxLC ={Ax1.1, …, Ax1.13}.

*Axiomas alternos* AxAL={Ax2.1, …,Ax2.4}. *Axiomas dobles* AxLD = AxLC ∩ AxAL.

Utilizando los operadores de verdad se infiere que:

Ax2.2 AxR: ¬X⊃~X ⇔ +X⊃X ⇔ X⊃⊗X.

Ax2.3 AxT: $\underline{X}$⊃+$\underline{X}$ ⇔ +X⊃++X ⇔ ¬X⊃+¬X ⇔ − ¬X⊃⊗X.

Definición 3. (*Teorema*). Para X∈FOR. Se dice que X es un *teorema* (X∈**TEO**) si y solamente si existe una *demostración de X*, es decir, X es la última fórmula de una sucesión finita de fórmulas, tales que cada una de ellas es un axioma o se infiere de dos fórmulas anteriores utilizando la regla de inferencia MP. El número de elementos de la sucesión se llama la *longitud de la prueba*.

En lo que sigue se utilizaran los siguientes resultados de la *Lógica Proposicional Clásica* LC, presentados en [7] y [8], y como en LD se tienen los correspondientes axiomas de CP y la regla de inferencia MP entonces, estos resultados valen en LD, donde X,Y,Z∈FOR.

I•c. *Introducción*: [(Z⊃X)•(Z⊃Y)]⊃[Z⊃(X•Y)].

I•. *Introducción de* •: X⊃(Y⊃(X•Y)).

E•c. *Eliminación*: [Z⊃(X•Y)]⊃[(Z⊃X)•(Z⊃Y)].

E•. *Eliminación de* •: (X•Y)⊃X, (X•Y)⊃Y.

SH. *Silogismo hipotético*: [(X⊃Y)•(Y⊃Z)]⊃(X⊃Z).

Exp. *Exportación*: [X⊃(Y⊃Z)]≡[(X•Y)⊃Z].

EQ. *Equivalencia*: (X≡Y)≡[(X⊃Y)•(Y⊃X)].

DN. *Doble negación*: X≡~~X.

DI. *Demostración indirecta*: X⊃(Y•~Y) ⇒ ~X.

N∪. *Negación de* ∪: ~(X∪Y)≡(~X•~Y).

N•. *Negación de* •: ~(X•Y)≡(~X∪~Y).

N⊃. *Negación de* ⊃: ~(X⊃Y)≡(X•~Y).

Tras. *Transposición*: (X⊃Y)≡(~Y⊃~X).

Imp. *Implicación*: (X⊃Y)≡(~X∪Y).

Id. *Principio de identidad :* X⊃X.

PR. *Principio de retorsión*: de (~X⊃X)≡X.

TD. *Teorema de deducción*: (X⇒Y) ⇒ (X⊃Y).

EQ. *Equivalencia*: (X≡Y)≡[(X•Y)∪(~X•~Y)].

SustEQ. F(X),X,Y∈FC y X≡Y ⇒ F(X)⊃F(Y).

Com•. *Conmutatividad*. (X•Y)≡(Y•X).

Aso•. *Asociatividad*. X•(Y•Z) ≡ (X•Y)•Z.

Idem•. *Idempotencia*. (X•X)≡X

Proposición 1 (*Construcción de verdades alternas*). Para X∈FOR.

X∈TEO entonces +X∈TEO.

Prueba: Supóngase que X∈TEO, se probará que +X∈TEO, haciendo inducción sobre la longitud de la demostración de X.

**PB**. Paso base. La longitud de la demostración de X es 1, es decir X es un axioma, pero si X es un axioma entonces, por Ax+, +X es axioma y por lo tanto, +X∈TEO.

**PI**. Paso de inducción. Como hipótesis inductiva **HI** se tiene que si la longitud de la demostración de Y es menor que L entonces +Y∈TEO. Supóngase que la demostración de X tiene longitud L mayor que 1. Se tiene entonces que X es un axioma o X es consecuencia de pasos anteriores utilizando la regla de inferencia MP. En el primer caso se procede como en el paso base. En el segundo caso se tienen, para alguna fórmula Z, demostraciones de Z⊃X y de Z, ambas demostraciones de longitud menor que L. Aplicando HI resulta que +(Z⊃X)∈TEO y +Z∈TEO. Como se tiene como axioma MP+: +(Z⊃X)⊃ (+Z⊃+X), aplicando dos veces la regla MP resulta que +X∈TEO. □

Proposición 2 (*Trasposiciones)*. Para X,Y,Z∈FOR. Si X⊃Y∈TEO entonces:

a. ~Y⊃~X∈TEO. b. ¬Y⊃¬X∈TEO. c. (X⊃~~Y), (~~X⊃Y)∈TEO.

d. ZX⊃ZY∈TEO e. ¬Y⊃¬(X•Y)∈TEO. f. (X•¬(X•Y))⊃¬Y∈TEO.

Prueba: Supóngase que X⊃Y∈TEO.

Parte a. Resulta de aplicar Tras en X⊃Y. Parte b. Por la parte a se infiere ~Y⊃~X∈TEO, por proposición 1 se afirma que +(~Y ⊃ ~X), por MP+ resulta +~Y ⊃ +~X, es decir, ¬~~Y ⊃ ¬~~X, lo cual por DN significa ¬Y⊃¬X∈TEO. Parte c. Resulta de aplicar parte a en parte a. Parte d. Supóngase Z•X, por E• se infieren Z y X, por MP resulta Y, por I• se afirma ZY, finalmente por TD se concluye Z•X⊃Z•Y∈TEO. □

Parte e. Por Ax1.7 se tiene que (X•Y)⊃Y∈TEO y en consecuencia por Tras también ~Y⊃~(X•Y)∈TEO, por la proposición 1 resulta +(~Y⊃~(X•Y))∈TEO, utilizando MP+ y MP se infiere +~Y⊃+~(X•Y)∈TEO, es decir,¬~~Y⊃¬~~(X•Y)∈TEO, por DN se concluye que ¬Y⊃¬(X•Y)∈TEO.

Parte b. Supóngase que X, y además supóngase que ¬(X•Y), por AxR y MP resulta ~(X•Y), por N• se infiere ~X∪~Y, y por DN y SD se deduce ~Y, utilizando TD resulta que ¬(X•Y)⊃~Y, por proposición 1 se obtiene +(¬(X•Y)⊃~Y), aplicando MP+ y MP deduce +¬(X•Y)⊃+~Y, es decir, +¬(X•Y)⊃¬Y, por AxT se tiene ¬(X•Y)⊃+¬(X•Y), aplicando SH se concluye ¬(X•Y)⊃¬Y. Finalmente, por TD se ha probado X⊃(¬(X•Y)⊃¬Y)∈TEO, lo cual por Exp equivale a (X•¬(X•Y))⊃¬Y∈TEO. □

# Principios de inducción estructural para LD

Definición 4. (*Número de negaciones*). Se define la función F de FOR en FOR de la siguiente manera ($Y, Y_0, Y_1, \ldots, Y_{n+1} \in FOR$):

$F_0(Y) = Y_0 \bullet Y$. $F_1(Y) = Y_1 \bullet \#F_0(Y) = Y_1 \bullet \#(Y_0 \wedge Y)$. $F_{n+1}(Y) = Y_{n+1} \bullet \#F_n(Y)$.

Donde $Z = \#Y$ significa que $Z = \neg Y$ o $Z = \sim Y$. $F_n(Y)$ significa: Y se encuentra en una fórmula rodeada por **n** negaciones, esto se denota como $\mathbf{Y^{Rn}}$.

**Notación:** Cuando sea conveniente y con el fin de simplificar la escritura, en lo que sigue, **XY** significa $X \bullet Y$, además por Aso•, XYZ significa X(YZ) o (XY)Z, también por Com• es indiferente escribir XY o YX.

Proposición 3. **PIEF-ID**. (*Principio de inducción estructural fuerte para iteración y desiteración*).

$(\forall n \in Z^+ \cap \{0\})[\underline{X} \in FA \text{ y } Y \in FOR \Rightarrow \underline{X} \bullet F_n(Y) \equiv \underline{X} \bullet F_n(\underline{X} \bullet Y)\ (\underline{X} \bullet Y^{Rn} \equiv \underline{X} \bullet (\underline{X} \bullet Y)^{Rn})]$

Prueba. Por inducción matemática sobre $n \in Z^+ \cap \{0\}$.

**PB.** Paso base. Sea n=0. Se sabe que $\underline{X}F_0(Y) = \underline{X}Y_0Y$. También se sabe que $\underline{X}F_0(Y\underline{X}) = \underline{X}Y_0\underline{X}Y$. La igualdad se satisface por las reglas E•, I• y Com•, junto con TD. Resulta entonces que $\underline{X} \bullet Y^{R0} \equiv \underline{X} \bullet (\underline{X} \bullet Y)^{R0}$.

Sea n=1. Se tiene que $\underline{X}F_1(Y) = \underline{X}Y_1\#F_0(Y) = \underline{X}Y_1\#\{Y_0Y\}$. Es decir, $\underline{X}F_1(Y) = \underline{X}Y_1\neg(Y_0Y)$ o $\underline{X}F_1(Y) = \underline{X}Y_1\sim(Y_0Y)$. También se sabe que $\underline{X}F_1(Y\underline{X}) = \underline{X}Y_1\#\{F_0(\underline{X}Y)\} = \underline{X}Y_1\#\{Y_0\underline{X}Y\}$, es decir, $\underline{X}F_1(Y\underline{X}) = \underline{X}Y_1\neg(Y_0\underline{X}Y)$ o $\underline{X}F_1(Y\underline{X}) = \underline{X}Y_1\sim(Y_0\underline{X}Y)$. De $\underline{X}Y_1\sim(Y_0Y)$, por E• se obtienen $\underline{X}Y_1$ y $\sim(Y_0Y)$, por N• resulta $\sim Y_0 \cup \sim Y$, por Ad se infiere $\sim Y_0 \cup \sim \underline{X} \cup \sim Y$, por N• resulta $\sim(Y_0\underline{X}Y)$, finalmente por I• se concluye $\underline{X}Y_1\sim(Y_0\underline{X}Y)$. Por TD se ha probado que $\underline{X}Y_1\sim(Y_0Y) \supset \underline{X}Y_1\sim(Y_0\underline{X}Y)$.

De $\underline{X}Y_1\sim(Y_0\underline{X}Y)$, por E• se obtienen $\underline{X}$, $Y_1$ y $\sim(Y_0\underline{X}Y)$, por N• resulta $\sim Y_0 \cup \sim \underline{X} \cup \sim Y$, aplicando SD se infiere $\sim Y_0 \cup \sim Y$, por N• resulta $\sim(Y_0Y)$, finalmente por I• se concluye $\underline{X}Y_1\sim(Y_0Y)$. Por TD se ha probado que $\underline{X}Y_1\sim(Y_0\underline{X}Y) \supset \underline{X}Y_1\sim(Y_0Y)$. De los dos últimos resultados por EQ se concluye que $\underline{X}Y_1\sim(Y_0\underline{X}Y) \equiv \underline{X}Y_1\sim(Y_0Y)$.

De $\underline{X}Y_1\neg(Y_0Y)$, por E• se obtienen $\underline{X}$, $Y_1$ y $\neg(Y_0Y)$, por la proposición 2e se infiere $\neg(Y_0\underline{X}Y)$, por I• se deduce $\underline{X}Y_1\neg(Y_0\underline{X}Y)$. Por TD se ha probado que $\underline{X}Y_1\neg(Y_0Y) \supset \underline{X}Y_1\neg(Y_0\underline{X}Y)$.

De $\underline{X}Y_1\neg(Y_0\underline{X}Y)$, por E• se obtienen $Y_1$ y $\underline{X}\neg(Y_0\underline{X}Y)$, utilizando la proposición 2f se deduce $\neg(Y_0Y)$, por I• resulta $\underline{X}Y_1\neg(Y_0Y)$, aplicando TD se concluye $\underline{X}Y_1\neg(Y_0\underline{X}Y) \supset \underline{X}Y_1\neg(Y_0Y)$, y como ya se probó la recíproca, por EQ se infiere que $\underline{X}Y_1\neg(Y_0\underline{X}Y) \equiv \underline{X}Y_1\neg(Y_0Y)$. Resulta entonces que $\underline{X} \bullet Y^{R1} \equiv \underline{X} \bullet (\underline{X} \bullet)^{R1}$.

**PI**. Paso inductivo. Hipótesis inductiva **HI**: $\underline{X}F_n(Y) \equiv \underline{X}F_n(\underline{X}Y)$, es decir, $\underline{X} \bullet Y^{Rn} \equiv \underline{X} \bullet (\underline{X} \bullet Y)^{Rn}$ n>1.

$\underline{X}F_{n+1}(Y) = \underline{X}Y_{n+1}\#F_n(Y)$, por el paso base equivale a $\underline{X}Y_{n+1}\#\{\underline{X}F_n(Y)\}$, por **HI** equivale a $\underline{X}Y_{n+1}\#\{\underline{X}F_n(\underline{X}Y)\}$, por el paso base equivale a $\underline{X}Y_{n+1}\#F_n(\underline{X}Y)$, es decir $\underline{X}F_{n+1}(\underline{X}Y)$. Se ha probado de esta manera que $\underline{X}F_{n+1}(Y) \equiv \underline{X}F_{n+1}(\underline{X}Y)$, es decir, $\underline{X} \bullet Y^{Rn+1} \equiv \underline{X} \bullet (\underline{X} \bullet Y)^{Rn+1}$

En consecuencia, por el principio de inducción matemática, se concluye PIEF-LD. □

Proposición 4. **PIEC-ID**. (*Principio de inducción estructural clásico para iteración y desiteración*).

$(\forall n\in Z^{+}\cap\{0\})[X,Y\in FOR \Rightarrow X\bullet G_n(Y) \equiv X\bullet G_n(X\bullet Y)$ y $X\bullet Y^{Rn} \equiv X\bullet(X\bullet Y)^{Rn}]$.

Donde $G_0(Y) = Y_0\bullet Y$ y $G_{n+1}(Y) = Y_{n+1}\bullet\sim G_n(Y)$. G es la restricción de F a solo negaciones clásicas.

**Observación**: en el PIEC-ID no se requiere que $X\in FA$, basta que $X\in FOR$.

Prueba. Basta observar que PIEC-ID es un caso particular de la proposición 3 PIEF-ID. □

Proposición 5. **PIE-RPI**. (*Principio de inducción estructural para regiones pares e impares*).

a. $(\forall n\in Z^{+}\cap\{0\})[X,Y\in FOR$ y $X\supset Y\in TEO \Rightarrow F_{2n}(X)\supset F_{2n}(Y)\in TEO$ $(X^{R2n} \supset Y^{R2n}\in TEO)]$

b. $(\forall n\in Z^{+}\cap\{0\})[X,Y\in FOR$ y $X\supset Y\in TEO \Rightarrow F_{2n+1}(Y)\supset F_{2n+1}(X)\in TEO$ $(Y^{R2n+1} \supset X^{R2n+1}\in TEO)]$

Prueba. Inducción matemática sobre $n\in Z^{+}\cap\{0\}$. Supóngase que $X\supset Y\in TEO$.

**PB.** Paso base. n=0. a. Por la proposición 2d se infiere $Y_0X\supset Y_0Y\in TEO$, es decir, $F_0(X)\supset F_0(Y)\in TEO$. Resulta entonces que $X^{R0}\supset Y^{R0}$.

b. Por la proposición 2 partes a, b, d se infiere $Y_1\{Y_0Y\}\supset Y_1\{Y_0X\}\in TEO$, es decir, $F_1(Y)\supset F_1(X)\in TEO$. Resulta entonces que $Y^{R1}\supset X^{R1}$.

**PI.** Paso inductivo. Parte a. Hipótesis inductiva **HI**: $F_{2n}(X)\supset F_{2n}(Y)\in TEO$. Por la proposición 2b se infiere $\neg(F_{2n}(Y))\supset \neg(F_{2n}(X))\in TEO$, por la proposición 2a resulta $\sim(F_{2n}(Y))\supset \sim(F_{2n}(X))\in TEO$, por lo que se tiene $\#F_{2n}(Y)\supset\#F_{2n}(X)\in TEO$, por la proposición 2d se concluye $Y_{2n+1}\#F_{2n}(Y)\supset Y_{2n+1}\#F_{2n}(X)\in TEO$, es decir, $F_{2n+1}(Y)\supset F_{2n+1}(X)\in TEO$, de nuevo la proposición 2 partes a, b, d implican $Y_{2n+2}\#F_{2n+1}(Y)\supset Y_{2n+2}\#F_{2n+1}(X)\in TEO$, lo cual significa $F_{2(n+1)}(X)\supset F_{2(n+1)}(Y)\in TEO$. Por lo tanto, $X^{R2(n+1)}\supset Y^{R2(n+1)}$.

Parte b. Hipótesis inductiva **HI**: $F_{2n+1}(Y)\supset F_{2n+1}(X)\in TEO$. Por la proposición 2 partes a, b, d se infiere $Y_{2n+2}\#F_{2n+1}(X)\supset Y_{2n+2}\#F_{2n+1}(Y)\in TEO$, es decir, $F_{2n+2}(X)\supset F_{2n+2}(Y)\in TEO$, de nuevo la proposición 2 partes a, b, d implican $Y_{2n+3}\#F_{2n+2}(Y)\supset Y_{2n+3}\#F_{2n+2}(X)\in TEO$, lo cual significa $F_{2(n+1)+1}(Y)\supset F_{2(n+1)+1}(X)\in TEO$. Resulta entonces que $Y^{R2(n+1)+1}\supset X^{R2(n+1)+1}$.

En consecuencia, por el principio de inducción matemática, se ha probado PIE-RPI. □

# Gráficos existenciales Gamma-LD

Se toma como punto de partida la *Hoja de aserción* **H** donde se dibujan los gráficos existenciales. El conjunto **Alfa-LC** de *gráficos existenciales clásicos*, se define de la siguiente manera:

1. $a\in FAC \Rightarrow a\in$ Alfa-LC.
2. $\lambda\in$ Alfa-LC ($\lambda$ es el gráfico vacío).
3. $X\in$ Alfa-LC $\Rightarrow$ **(**X**)**$\in$ Alfa-LC.
4. $X,Y\in$ Alfa-LC $\Rightarrow XY\in$ Alfa-LC.
5. Sólo 1 a 4 determinan Alfa-LC.

El conjunto **Alfa-LI** de *gráficos existenciales intuicionistas*, se define de la siguiente manera:

1. a∈FAA ⇒ a∈Alfa-LI.
2. λ∈Alfa-LI (λ es el gráfico vacío)
3. X∈Alfa-LI ⇒ **[**X**]**∈Alfa-LI.
4. X,Y∈Alfa-LI ⇒ XY, **[**X **(**Y**)]**, **[(**X**)** **(**Y**)]**∈Alfa-LI.
5. Sólo 1 a 4 determinan Alfa-LI.

El conjunto **Gamma-LD** de ***gráficos existenciales LD*** se define de la siguiente manera:

1. X∈Alfa-LC ⇒ X∈Gamma-LD.
2. λ∈Gamma-LD (λ es el gráfico vacío).
3. X∈Alfa-LI ⇒ X∈Gamma-LD.
4. X∈Gamma-LD ⇒ **(**X**)**,**[**X**]**∈Gamma-LD.
5. X,Y∈Gamma-LD ⇒ XY,**[**X **(**Y**)]**,**[(**X**)** **(**Y**)]**∈Gamma-LD.
6. Sólo 1 a 5 determinan Gamma-LD.

El conjunto **GA** de **gráficos existenciales alternos** se define de la siguiente manera:

X∈GA ⇔ X=**[**Y**]** con Y∈Gamma-LD o X=a donde a∈FAA

Notación. Si X∈GAE entonces **X** indica que X∈GA.

# Reglas de transformación y Validez

Cuando se tiene el gráfico existencial **(**X**)**, se dice X está rodeado por un *corte clásico* y cuando se tiene el gráfico existencial **[**X**]**, se dice X está rodeado por un *corte alterno*. Se dice que un gráfico existencial X se encuentra en una *región par*, denotado **X∈RP**, si X se encuentra rodeado por un número par de cortes (clásicos y/o alternos). X se encuentra en una *región impar*, denotado **X∈RI**, si X se encuentra rodeado por un número impar de cortes (clásicos y/o alternos). **X∈R$_k$P** significa que X se encuentra en la *región k*, la cual es par. **X∈R$_k$I** significa que X se encuentra en la región k, la cual es impar. X se encuentra en una *región clásica*, denotado **X∈RC**, si X no se encuentra rodeado por cortes alternos. X se encuentra en una *región alterna*, denotado **X∈RA**, si X se encuentra rodeado por al menos un corte alterno.

Notación. $\mathbf{X} \overset{\boldsymbol{Re}}{\Rightarrow} \mathbf{Z}$ significa que la *regla de transformación Re*, aplicada a un gráfico X permite inferir un nuevo gráfico Z. Si también Z se transforma en X mediante *Re* entonces se escribe $\mathbf{X} \overset{\boldsymbol{Re}}{\Leftrightarrow} \mathbf{Z}$.

**X>>Z** significa que X se transforma en Z utilizando un número finito de reglas de transformación. Para X∈Gamma-LD, se dice que X es *válido*, denotado X∈**GEV**, si λ>>X.

Definición 5. **(***Definición de λ*).

**Def-λ**. Definición de λ. $\lambda \overset{\mathbf{Def-\lambda}}{\Longleftrightarrow}$ ‘ ’ (gráfico vacío).

El conjunto **RTRA** de *reglas primitivas de transformación* consta de (donde X,Y,Z∈Gamma-LD):

1. **Rλ**: Regla lambda. λ∈GEV

2. **B:** Borrado. $XY \in RP \Rightarrow XY \overset{\mathbf{B}}{\Rightarrow} X$

$XY \in RP \Rightarrow XY \overset{\mathbf{B}}{\Rightarrow} Y$

3. **E:** Escritura. $X \in RI \Rightarrow X \overset{\mathbf{E}}{\Rightarrow} XY$

$X \in RI \Rightarrow X \overset{\mathbf{E}}{\Rightarrow} YX$

4. **DCC**: Doble corte clásico. $X \overset{\mathbf{DCC}}{\Leftrightarrow} ((X))$

5. **CC**: Cambio de corte. $[X] \in RP \Rightarrow [X] \overset{\mathbf{CC}}{\Rightarrow} (X)$

$(X) \in RI \Rightarrow (X) \overset{\mathbf{CC}}{\Rightarrow} [X]$

6. **DCMGEV**: Doble corte mixto para GEV.

$X \in GEV \Rightarrow X \overset{\mathbf{DCMGEV}}{\Longleftrightarrow} [(X)]$

7. **DCMF**: Doble corte mixto fuerte. $\underline{X} \in RP \Rightarrow \underline{X} \overset{\mathbf{DCMF}}{\Longrightarrow} [(\underline{X})]$

$\underline{X} \in RI \Rightarrow [(\underline{X})] \overset{\mathbf{DCMF}}{\Longrightarrow} \underline{X}$

8. **I**: Iteración. $X \overset{\mathbf{I}}{\Rightarrow} XX$

9. **D**: Desiteración. $XX \overset{\mathbf{D}}{\Rightarrow} X$

10. **IC**: Iteración clásica. $X,Y \in RC \Rightarrow X(\ldots(Y)\ldots) \overset{\mathbf{IC}}{\Longrightarrow} X(\ldots(XY)\ldots)$

11. **DC**: Desiteración clásica. $X,Y \in RC \Rightarrow X(\ldots(XY)\ldots) \overset{\mathbf{DC}}{\Longrightarrow} X(\ldots(Y)\ldots)$

12. **IF**: Iteración fuerte. $\underline{X}\{\ldots\{Y\}\ldots\} \overset{\mathbf{IF}}{\Longrightarrow} \underline{X}\{\ldots\{\underline{X}Y\}\ldots\}$

Donde {Z} significa (Z) o [Z].

13. **DF**: Desiteración fuerte. $\underline{X}\{\ldots\{\underline{X}Y\}\ldots\} \overset{\mathbf{DF}}{\Longrightarrow} \underline{X}\{\ldots\{Y\}\ldots\}$

14. Sólo las reglas 1 a 13 determinan el conjunto **RTRA**.

Notación. Las reglas de la forma $X \overset{\mathbf{R}}{\Leftrightarrow} Y$ indican que $X \overset{\mathbf{R}}{\Longrightarrow} Y$ y $Y \overset{\mathbf{R}}{\Longrightarrow} X$ valen tanto en regiones pares como en regiones impares.

Observaciones.

A. Formalmente, las reglas IC y DC se expresan así: **IC**: $XG_n(Y) \overset{\mathbf{IC}}{\Longrightarrow} XG_n(XY)$. **DC**: $XG_n(XY) \overset{\mathbf{DC}}{\Longrightarrow} XG_n(Y)$. Donde $G_0(Y)=Y_0Y$, $G_1(Y) = Y_1\mathbf{(}G_0(Y)\mathbf{)} = Y_1\mathbf{(}Y_0Y\mathbf{)}$, $G_{n+1}(Y) = Y_{n+1}\mathbf{(}G_n(Y)\mathbf{)}$. $G_n(Y)$ significa que, Y se encuentra en una región rodeada por **n** cortes clásicos, esto se denota como $\mathbf{Y^{Rn}}$.

B. Formalmente, las reglas IF y DF se expresan así: **IF**: $XH_n(Y) \xRightarrow{\textbf{IF}} XH_n(XY)$. **DF**: $XH_n(XY) \xRightarrow{\textbf{DF}} XH_n(Y)$. Donde $H_0(Y) = Y_0Y$, $H_1(Y) = Y_1\mathbf{[}H_0(Y)\mathbf{]} = Y_1\mathbf{[}Y_0Y\mathbf{]}$, $H_{n+1}(Y) = Y_{n+1}\mathbf{[}H_n(Y)\mathbf{]}$. $H_n(Y)$ significa que, Y se encuentra en una región rodeada por **n** cortes intuicionistas, se denota como $\mathbf{Y^{Rn}}$.

Proposición 6. (*Reglas de transformación derivadas*).

1. **DCCλ**: Doble corte clásico λ. — $((\lambda))\in GEV$, $((\ ))\in GEV$.
2. **DCM**: Doble corte mixto. — $X\in RP \Rightarrow [(X)] \xRightarrow{\textbf{DCM}} X$.
   $X\in RI \Rightarrow X \xRightarrow{\textbf{DCM}} [(X)]$.
3. **DCMλ**: Doble corte mixto λ. — $[(\lambda)]\in GEV$, $[(\ )]\in GEV$.
4. **CCE**: Cambio de corte especifico. — $[X\ (Y)]\in RP \Rightarrow [X\ (Y)] \xRightarrow{\textbf{CCE}} [X\ [Y]]$.
   $[X\ [Y]]\in RI \Rightarrow [X\ [Y]] \xRightarrow{\textbf{CCE}} [X\ (Y)]$.
5. **DCMF.1**: Doble corte mixto fuerte. — $[(\underline{X})] \xLeftrightarrow{\textbf{DCMF.1}} \underline{X}$.
6. **TCM**: Triple corte mixto. — $[([X])] \xLeftrightarrow{\textbf{TCM}} [X]$.
7. **DCAF**: Doble corte alterno fuerte. — $\underline{X}\in RP \Rightarrow \underline{X} \xRightarrow{\textbf{DCAF}} [[\underline{X}]]$.
   $\underline{X}\in RI \Rightarrow [[\underline{X}]] \xRightarrow{\textbf{DCAF}} \underline{X}$.
8. **TCA**: Triple corte alterno. — $[X]\in RP \Rightarrow [X] \xRightarrow{\textbf{TCA}} [[[X]]]$.
   $[X]\in RI \Rightarrow [[[X]]] \xRightarrow{\textbf{TCA}} [X]$.
9. **TCAF**: Triple corte alterno fuerte. — $[\underline{X}]\in RP \Rightarrow [[[\underline{X}]]] \xRightarrow{\textbf{TCAF}} [\underline{X}]$
   $[\underline{X}]\in RI \Rightarrow [\underline{X}] \xRightarrow{\textbf{TCAF}} [[[\underline{X}]]]$
10. **TCAF.1**: Triple corte alterno fuerte. — $[[[\underline{X}]]] \xLeftrightarrow{\textbf{TCAF.1}} [\underline{X}]$.
11. **CCA**: Cuádruple corte alterno. — $[[[[X]]]] \xLeftrightarrow{\textbf{CCA}} [[X]]$.
12. **ID**: Iteración. Desiteración. — $X \xLeftrightarrow{\textbf{ID}} XX$.

Prueba 1: $\lambda \xLeftrightarrow{\textbf{DCC}} ((\lambda))$, y por Rλ se tiene que $\lambda\in GEV$, por lo que $((\lambda))\in GEV$, por B resulta $((\ ))\in GEV$.

Prueba 2: $X\in RP \Rightarrow [(X)]\in RP \Rightarrow [(X)] \xRightarrow{\textbf{CC}} ((X)) \xRightarrow{\textbf{DCC}} X$. $X\in RI \Rightarrow ((X))\in RI \Rightarrow X \xRightarrow{\textbf{DCC}} ((X)) \xRightarrow{\textbf{CC}} [(X)]$.

Prueba 3: Por Rλ se tiene $\lambda\in GEV$, aplicando DCMGEV resulta $\mathbf{[(}\lambda\mathbf{)]}$, por B se infiere, $\mathbf{[(\ )]}\in GEV$.

Prueba 4: $[X\ (Y)]\in RP \Rightarrow (Y)\in RI \Rightarrow [X\ (Y)] \xRightarrow{\textbf{CC}} [X\ [Y]]$. $[X\ [Y]]\in RI \Rightarrow [Y]\in RP \Rightarrow [X\ [Y]] \xRightarrow{\textbf{CC}} [X\ (Y)]$

Prueba 5: Se tienen $\underline{X}\in RP \Rightarrow \underline{X} \xRightarrow{\mathbf{DCMF}} [(\underline{X})]$ y $\underline{X}\in RP \Rightarrow [(\underline{X})] \xRightarrow{\mathbf{DCM}} \underline{X}$, también se tienen

$\underline{X}\in RI \Rightarrow [(\underline{X})] \xRightarrow{\mathbf{DCMF}} \underline{X}$ y $\underline{X}\in RI \Rightarrow \underline{X} \xRightarrow{\mathbf{DCM}} [(\underline{X})]$. Por lo que, $[(\underline{X})] \xLeftrightarrow{\mathbf{DCMF.1}} \underline{X}$.

Prueba 6: En particular se tiene $[([X])] \xLeftrightarrow{\mathbf{DCMF.1}} [X]$.

Prueba 7: $\underline{X}\in RP \Rightarrow \underline{X} \xRightarrow{\mathbf{DCMF}} [(\underline{X})] \xRightarrow{\mathbf{CC}} [[\underline{X}]]$, además, $\underline{X}\in RI \Rightarrow [\underline{X}]\in RP \Rightarrow [[\underline{X}]] \xRightarrow{\mathbf{CC}} [(\underline{X})] \xRightarrow{\mathbf{DCMF}} \underline{X}$.

Prueba 8: Caso particular de DCAF.

Prueba 9: $[\underline{X}]\in RP \Rightarrow [[[\underline{X}]]] \xRightarrow{\mathbf{CC}} [[(\underline{X})]] \xRightarrow{\mathbf{DCMF}} [\underline{X}]$. Además, $[X]\in RI \Rightarrow [[[X]]] \xRightarrow{\mathbf{TCA}} [X]$.

Prueba 10: $[\underline{X}]\in RP \Rightarrow [[[\underline{X}]]] \xRightarrow{\mathbf{TCAF}} [\underline{X}]$ y $[\underline{X}]\in RP \Rightarrow [\underline{X}] \xRightarrow{\mathbf{TCA}} [[[\underline{X}]]]$. Luego, $[\underline{X}]\in RP \Rightarrow [\underline{X}] \Leftrightarrow [[[\underline{X}]]]$.

$[\underline{X}]\in RI \Rightarrow [\underline{X}] \xRightarrow{\mathbf{TCAF}} [[[\underline{X}]]]$ y $[\underline{X}]\in RI \Rightarrow [[[\underline{X}]]] \xRightarrow{\mathbf{TCA}} [\underline{X}]$. Luego, $[\underline{X}]\in RI \Rightarrow [[[\underline{X}]]] \Leftrightarrow [\underline{X}]$. Por lo tanto, $[\underline{X}] \Leftrightarrow [[[\underline{X}]]]$.

Prueba 11: Caso particular de TCAF.1.

Prueba 12: $X \xRightarrow{\mathbf{I}} XX$. Además, $XX \xRightarrow{\mathbf{D}} X$. Por lo tanto, $X \Leftrightarrow XX$. □

Proposición 7. (*Reversión de las reglas de transformación*). Para $X,Y\in$Gamma-LD.

a. $(\forall R\in RTRA)[X\in RP$ y $X\xRightarrow{\mathbf{R}}Y](\exists R'\in RTRA)[Y\in RI$ y $Y\xRightarrow{R'}X]$

b. $(\forall R\in RTRA)\{[X\in R_kP$ y $X\xRightarrow{\mathbf{R}}Y \Rightarrow Y\in R_kP]$ y $[X\in R_kI$ y $X\xRightarrow{\mathbf{R}}Y \Rightarrow Y\in R_kI]\}$

Prueba: por simple inspección de las reglas primitivas 1 a 14. □

Proposición 8. (*Principio de contraposición*). Para $X_0,X_n\in$Gamma-LD.

$X_0>>X_n \Rightarrow \{(X_0\in RP \Rightarrow X_0>>X_n)$ y $(X_n\in RI \Rightarrow X_n>>X_0)\}$

Prueba: Si $X_0>>X_n$ entonces existen $R_1, \ldots, R_n\in RTRA$, y existen $X_1, \ldots, X_{n-1}\in$Gamma-LD, tales que $X_0R_1X_1R_2X_2\ldots X_{n-1}R_nX_n$, y se dice que la *longitud* de la transformación de $X_0>>X_n$ es **n**.

La prueba se realiza por inducción sobre la longitud de la transformación.

**PB**. Paso base. n=1. Significa que solo se aplicó una de las reglas primitivas, y como $X_0$ está en una región de paridad 0, entonces R tiene que ser de la forma $X_0\in RP \Rightarrow X_0\xRightarrow{\mathbf{R}}X_1$ . Por la proposición 7a se infiere que $X_1\in RI \Rightarrow X_1\xRightarrow{R'}X_0$.

**PI**. Paso inductivo. Hipótesis inductiva **HI**: $(\forall n>1)[X_0>>X_n \Rightarrow \{(X_0\in RP \Rightarrow X_0>>X_n)$ y $(X_n\in RI \Rightarrow X_n>>X_0)\}]$. Si $X_0>>X_{n+1}$ entonces $X_0R_1X_1R_2X_2\ldots X_{n-1}R_nX_nR_{n+1}X_{n+1}$, es decir, $X_0R_1X_1R_2X_2\ldots X_{n-1}R_nX_n$ y $X_nR_{n+1}X_{n+1}$, por lo que, $X_0>>X_n$ y $X_nR_{n+1}X_{n+1}$. Aplicando **HI** y la proposición 7b resultan $(X_0\in RP \Rightarrow X_0>>X_n$ y $X_n\in RP)$, $(X_n\in RP \Rightarrow X_nR_{n+1}X_{n+1}$ y $X_{n+1}\in RP)$, $(X_n\in RI \Rightarrow X_n>>X_0$ y $X_0\in RI)$ y $(X_{n+1}\in RI \Rightarrow X_{n+1}R'_{n+1}X_n$ y $X_n\in RI)$. Por lo que, $(X_0\in RP \Rightarrow X_0>>X_{n+1})$ y $(X_{n+1}\in RI \Rightarrow X_{n+1}>>X_0)$.

Por el principio de inducción matemática se concluye la veracidad del principio de contraposición. □

Proposición 9. TDG y TDGF. (*Teoremas de deducción gráfico y fuerte en Gamma-LD*).

Para $X,Y \in$ Gamma-LD. a. TDG. $X \gg Y \Rightarrow (X\ (Y))$ b. TDGF. $\underline{X} \gg \underline{Y} \Rightarrow [\underline{X}\ (\underline{Y})]$

Prueba. Supóngase que $X \gg Y$.

Parte a. $\lambda \xRightarrow{\mathbf{DCC\lambda}} ((\ )) \overset{\mathbf{E}}{\Rightarrow} (X\ (\ )) \overset{\mathbf{IC}}{\Rightarrow} (X\ (X)) \xRightarrow{\mathbf{X\gg Y}\text{ y proposición 8}} (X\ (Y))$. Por lo tanto, $X \gg Y \Rightarrow (X\ (Y))$.

Parte b. $\lambda \xRightarrow{\mathbf{DCM\lambda}} [(\ )] \overset{\mathbf{E}}{\Rightarrow} [\underline{X}\ (\ )] \overset{\mathbf{IF}}{\Rightarrow} [\underline{X}\ (\underline{X})] \xRightarrow{\mathbf{X\gg Y}\text{ y proposición 8}} [\underline{X}\ (\underline{Y})]$. Por lo tanto, $\underline{X} \gg \underline{Y} \Rightarrow [\underline{X}\ (\underline{Y})]$. □

Proposición 10. TDIG. (*Teorema de demostración indirecta en Gamma-LD*). Para $X \in$ Gamma-LD.

a. $X \gg (\ ) \Rightarrow (X)$ b. $(X) \gg (\ ) \Rightarrow X$ c. $\underline{X} \gg [\ ] \Rightarrow [\underline{X}]$

Prueba. Parte a. Supóngase que $X \gg (\ )$. Por TDG (proposición 9) se infiere $\mathbf{(}X\ \mathbf{((\ ))\ )} \xRightarrow{\mathbf{DCC}} \mathbf{(}X\mathbf{)}$.

Parte b. Supóngase que $(X) \gg (\ )$. Por la parte a se deduce $\mathbf{((}X\mathbf{))} \xRightarrow{\mathbf{DCC}} X$. □

Parte c. Supóngase que $\underline{X} \gg [\ ]$. Por TDGF (proposición 9) se infiere $[\underline{X}\ ([\ ])] \overset{\mathbf{CC}}{\Rightarrow} [\underline{X}\ [[\ ]]] \xRightarrow{\mathbf{DCAF}} [\underline{X}]$. □

# Traducción de FOR a Gamma-LD

Definición 6. (*Traducción de fórmulas a gráficos existenciales*). La función **( )’** de FOR en Gamma-LD, se define de la siguiente manera (donde $a \in$ **FAT** y X,Y):

1. a’ = a
2. {~X}’ = (X’)
3. {X•Y}’ = X’Y’
4. {X⊃Y}’ = (X’(Y’))
5. {X≡Y}’ = {X⊃Y}’ {Y⊃X}’ = (X’(Y’)) (Y’(X’))
6. {X∪Y}’ = ((X’) (Y’))
7. Ax1.1’ = λ
8. {¬X}’ = [X’]
9. {+X}’ = [(X’)]
10. {X→Y}’ = [X’ (Y’)]
11. {X∨Y}’ = [(X’) (Y’)]
12. {X∧Y}’ = [(X’ Y’)]
13. {X↔Y}’ = [X’ (Y’)] [Y’ (X’)]

Proposición 11. (*Validez en LD de las reglas primitivas de transformación en LD*). Para $X,Y \in$ FOR.

$(\forall R \in \text{RTRA})(\exists T \in \text{TEO})[X' \overset{\mathbf{R}}{\Rightarrow} Y' \Rightarrow T = X \supset Y]$

Prueba. Regla 1. Rλ: Regla lambda. $\lambda \in$ GEV.

$\lambda$=Ax1.1’ y Ax1.1$\in$TEO.

Regla 2. B: Borrado. $X'Y' \in \text{RP} \Rightarrow (X'Y' \overset{\mathbf{B}}{\Rightarrow} X')$ y $(X'Y' \overset{\mathbf{B}}{\Rightarrow} Y')$.

Regla 3. E: Escritura. $X'Y' \in \text{RI} \Rightarrow (X' \overset{\mathbf{E}}{\Rightarrow} X'\ Y')$ y $(Y' \overset{\mathbf{E}}{\Rightarrow} X'Y')$.

Por Ax1.6 y Ax1.7 se tienen $(X\bullet Y)\supset X\in TEO$ y $(X\bullet Y)\supset Y\in TEO$, por la proposición 5 se concluyen que $(XY)^{R2n}\supset Y^{R2n}\in TEO$ y $(XY)^{R2n}\supset X^{R2n}\in TEO$, $Y^{R2n+1}\supset(XY)^{R2n+1}\in TEO$ y $X^{R2n+1}\supset(XY)^{R2n+1}\in TEO$.

Regla 4. DCC: Doble corte clásico. $X' \overset{\mathbf{DCC}}{\Longleftrightarrow} ((X'))$.

Por DN se tiene $\sim\sim X\equiv X\in TEO$, aplicando la proposición 5 se concluye que $(\sim\sim X)^{R2n}\equiv X^{R2n}\in TEO$ y $(\sim\sim X)^{R2n+1}\equiv X^{R2n+1}\in TEO$, por lo que, $(\sim\sim X)^{Rn}\equiv X^{Rn}\in TEO$.

Regla 5. CC: Cambio de corte. $[X']\in RP \Rightarrow [X'] \overset{\mathbf{CC}}{\Rightarrow} (X')$ y $(X')\in RI \Rightarrow (X') \overset{\mathbf{CC}}{\Longrightarrow} [X']$.

Por AxR se tiene $\neg X\supset\sim X\in TEO$, por la proposición 5 se concluye que $(\neg X)^{R2n}\supset(\sim X)^{R2n}\in TEO$ y $(\sim X)^{R2n+1}\supset(\neg X)^{R2n+1}\in TEO$.

Regla 6. DCMGEV: Doble corte mixto para GEV. $X'\in GEV \Rightarrow X' \overset{\mathbf{DCMGEV}}{\Longleftrightarrow} [(X')]$.

Por la proposición 1 se sigue que $X\in TEO \Rightarrow X\supset +X\in TEO$, por AxR se tiene que $+X\supset X\in TEO$, resultando que $X\in TEO \Rightarrow X\equiv +X\in TEO$.

Regla 7. DCMF: Doble corte mixto fuerte. $\underline{X}'\in RP \Rightarrow \underline{X}' \overset{\mathbf{DCMF}}{\Longrightarrow} [(\underline{X}')]$ y $\underline{X}'\in RI \Rightarrow [(\underline{X}')] \overset{\mathbf{DCMF}}{\Longrightarrow} \underline{X}'$.

Por AxT se tiene $\underline{X}\supset +\underline{X}\in TEO$, aplicando la proposición 5 se concluye que $(\underline{X})^{R2n}\supset(+\underline{X})^{R2n}\in TEO$ y $(+\underline{X})^{R2n+1}\supset(\underline{X})^{R2n+1}\in TEO$.

Regla 8. I: Iteración. $X' \overset{\mathbf{I}}{\Rightarrow} X'X'$.

Regla 9. D: Desiteración. y $X'X' \overset{\mathbf{D}}{\Rightarrow} X'$.

Por Idem$\bullet$ se tiene $(X\bullet X)\equiv X$ aplicando la proposición 5 se concluye que $(XX)^{R2n}\equiv X^{R2n}\in TEO$ y $(XX)^{R2n+1}\equiv X^{R2n+1}\in TEO$.

Regla 10. IC: Iteración clásica. $X',Y'\in RC \Rightarrow X'(\ldots(Y')\ldots) \overset{\mathbf{IC}}{\Longrightarrow} X'(\ldots(X'Y')\ldots)$.

Regla 11. DC: Desiteración clásica. $X',Y'\in RC \Rightarrow X'(\ldots(X'Y')\ldots) \overset{\mathbf{DC}}{\Longrightarrow} X'(\ldots(Y')\ldots)$.

Por la proposición 4 se tiene que, para $n\in Z^{+}\cap\{0\}$, $X\bullet Y^{Rn} \equiv X\bullet(X\bullet Y)^{Rn}\in TEO$.

Regla 12. IF: Iteración fuerte. $\underline{X}'\{\ldots\{Y'\}\ldots\} \overset{\mathbf{IF}}{\Longrightarrow} \underline{X}'\{\ldots\{\underline{X}'Y'\}\ldots\}$.

Regla 13. DF: Desiteración fuerte. $\underline{X}'\{\ldots\{\underline{X}Y\}\ldots\} \overset{\mathbf{DF}}{\Longrightarrow} \underline{X}\{\ldots\{Y\}\ldots\}$, donde $\{Z\}$ significa $(Z)$ o $[Z]$.

Por la proposición 3 se tiene que, para $n\in Z^{+}\cap\{0\}$, $\underline{X}\bullet Y^{Rn} \equiv \underline{X}\bullet(\underline{X}\bullet Y)^{Rn}$. □

Proposición 12. (*Los gráficos válidos en Gamma-LD son teoremas en LD*).

$(\forall Z\in FOR)(Z'\in GEV \Rightarrow Z\in TEO)$

Prueba: Si $\lambda >> X_n$’ entonces existen $R_1, \ldots, R_n \in$RTRA, y existen $\lambda$, $X_2$’, …, $X_{n-1}$’$\in$Gamma-LD, tales que $\lambda R_1 X_1$’$R_2 X_2$’…$X_{n-1}$’$R_n X_n$’, y se dice que la *longitud* de la transformación de $\lambda >> X_n$’ es **n**.

La prueba se realiza por inducción sobre la longitud *n* de la transformación.

**PB**. Paso base. *n*=1. Significa que solo se aplicó una de las reglas primitivas, y $R_1$ tiene que ser de la forma $\lambda \overset{\mathbf{R1}}{\Longrightarrow} X_1$’ . Por la proposición 11 se infiere que Ax1.1$\supset X_1 \in$TEO, luego $X_1 \in$TEO.

**PI**. Paso inductivo. Hipótesis inductiva **HI**: $(\forall n>1)[\lambda >> X_n$’ $\Rightarrow X_n \in$TEO. Si $\lambda >> X_{n+1}$’ entonces $\lambda R_1 X_1$’$R_2 X_2$’…$X_{n-1}$’$R_n X_n$’$R_{n+1} X_{n+1}$’, es decir, $\lambda R_1 X_1$’$R_2 X_2$’…$X_{n-1}$’$R_n X_n$’ y $X_n$’$R_{n+1} X_{n+1}$’, por lo que, $\lambda >> X_n$’ y $X_n$’$R_{n+1} X_{n+1}$’. Aplicando **HI** y TDG proposición 11 se infieren $X_n \in$TEO y $X_n \supset X_{n+1} \in$TEO, aplicando MP se concluye que $X_{n+1} \in$TEO. Por el principio de inducción matemática se tiene que $(\forall Z \in$FOR$)(Z$’$\in$GEV $\Rightarrow Z \in$TEO). □

Proposición 13. (*AxLC válidos en Gamma-LD*).

$(\forall X \in$AxLC$)(X$’$\in$GEV)

Prueba. Ax1.1. $\lambda \overset{\mathbf{DCC\lambda}}{\Longrightarrow} ((\lambda)) \overset{\mathbf{DCC\lambda}}{\Longrightarrow} ((\ ((\ ))\ )) \overset{\mathbf{E}}{\Rightarrow} (X’(\ (Y’(\ ))\ )) \overset{\mathbf{IC}}{\Rightarrow} (X’(\ (Y’(X’))\ ))$. Por lo tanto, (Ax1.1)’$\in$GEV.

Ax1.2. $(X’((Y’(Z’))))\ (X’(Y’))\ X’ \overset{\mathbf{DCC}}{\Longrightarrow} (X’Y’(Z’))\ (X’(Y’))\ X’ \overset{\mathbf{DC}}{\Longrightarrow} (Y’(Z’))\ ((Y’))\ X’ \overset{\mathbf{DCC}}{\Longrightarrow} (Y’(Z’))\ Y’\ X’ \overset{\mathbf{D}}{\Rightarrow} ((Z’))\ Y’\ X’ \overset{\mathbf{DCC}}{\Longrightarrow} Z’\ Y’\ X’ \overset{\mathbf{B}}{\Rightarrow} Z’$. Aplicando TDG (proposición 9) 3 veces resulta que (Ax1.2)’$\in$GEV.

Ax1.3 y Ax1.4. $\lambda \overset{\mathbf{DCC\lambda}}{\Longrightarrow} ((\lambda)) \overset{\mathbf{DCC\lambda}}{\Longrightarrow} ((\ ((\ ))\ )) \overset{\mathbf{E}}{\Rightarrow} (X’(\ ((Y’)\ (\ ))\ )) \overset{\mathbf{IC}}{\Rightarrow} (X’(\ ((Y’)\ (X’))\ ))$. Por lo tanto,

(Ax1.3)’$\in$GEV y (Ax1.4)’$\in$GEV.

Ax1.5. $(X’(Z’))\ (Y’(Z’))\ ((X’)(Y’))\ (Z’) \overset{\mathbf{DC}}{\Longrightarrow} (X’)\ (Y’)\ ((X’)(Y’))\ (Z’) \overset{\mathbf{DC}}{\Longrightarrow} (X’)\ (Y’)\ (\ )\ (Z’) \overset{\mathbf{B}}{\Rightarrow} (\ )$. Por TDIG (proposición 10) se infiere que $(X’(Z’))\ (Y’(Z’))\ ((X’)(Y’)) \gg Z’$. Aplicando TDG 3 veces resulta que (Ax1.5)’$\in$GEV.

Ax1.6 y Ax1.7. $\lambda \overset{\mathbf{DCC\lambda}}{\Longrightarrow} ((\ )) \overset{\mathbf{E}}{\Rightarrow} (X’Y’(\ )) \overset{\mathbf{IC}}{\Rightarrow} (X’Y’(X’))$. Por lo tanto, (Ax1.6)’$\in$GEV y (Ax1.7)’$\in$GEV.

Ax1.8. $(X’(Y’))\ (X’(Z’))\ X’ \overset{\mathbf{DC}}{\Longrightarrow} ((Y’))\ ((Z’))\ X’ \overset{\mathbf{DCC}}{\Longrightarrow} Y’\ Z’\ X’$. Aplicando TDG 2 veces resulta que (Ax1.8)’$\in$GEV.

Ax1.9. La misma prueba de Ax1.3. Por lo tanto, (Ax1.9)’$\in$GEV.

Ax1.10. $\lambda \overset{\mathbf{DCC\lambda}}{\Longrightarrow} ((\ )) \overset{\mathbf{E}}{\Rightarrow} (X’\ (\ )) \overset{\mathbf{I}}{\Rightarrow} (X’\ (X’)) \overset{\mathbf{DCC}}{\Longrightarrow} (((X’))\ (X’))$. Por lo tanto, (Ax1.10)’$\in$GEV.

Ax1.11 y Ax1.12. $(X \equiv Y)’ \overset{\text{Definición}}{\Longrightarrow} (X \supset Y)’\ (Y \supset X)’ \overset{\mathbf{B}}{\Rightarrow} (Y \supset X)$. Por lo tanto, (Ax1.11)’$\in$GEV y (Ax1.12)’$\in$GEV.

Ax1.13. $(X \supset Y)’\ (Y \supset X)’ \overset{\text{Definición}}{\Longrightarrow} (X \equiv Y)’$. Aplicando TDG 2 veces resulta que (Ax1.13)’$\in$GEV. □

Proposición 14. (*AxAL válidos en Gamma-LD*).

(∀X∈AxAL)(X’∈GEV)

Prueba. Ax2.1. [((X’(Y’)))] [(X’)]([(Y’)]) $\overset{\mathbf{DCC}}{\Longrightarrow}$ [X’(Y’)] [(X’)]([(Y’)]) $\overset{\mathbf{DCM}}{\Longrightarrow}$ [[(X’)](Y’)] [(X’)]([(Y’)]) $\overset{\mathbf{DF}}{\Longrightarrow}$ [(Y’)] [(X’)] ([(Y’)]) $\overset{\mathbf{B}}{\Rightarrow}$ [(Y’)] ([(Y’)]) $\overset{\mathbf{D}}{\Rightarrow}$ [(Y’)] ( ) $\overset{\mathbf{B}}{\Rightarrow}$ ( ). Por TDIG (proposición 10) se infiere que [((X’(Y’)))] [(X’)] >> [(Y’)]. Aplicando TDG (proposición 9) 2 veces resulta que (Ax2.1)’∈GEV.

Ax2.2. [X’] $\overset{\mathbf{CC}}{\Rightarrow}$ (X’). Aplicando TDG resulta que (Ax2.2)’∈GEV.

Ax2.3. [X’] $\overset{\mathbf{TCM}}{\Longrightarrow}$ [([X’])]. Aplicando TDG resulta que (Ax2.3)’∈GEV.

Ax2.4. Se acaba de probar que si X∈{Ax1.1, …, Ax1.13, Ax2.1, …, Ax2.3} entonces X’∈GEV $\overset{\mathbf{DCMGEV}}{\Longrightarrow}$ [(X’)]. Por lo tanto, (Ax2.4)’∈GEV. □

Proposición 15. (*MP)’ preserva validez en Gamma-LD*).

(∀X,Y∈FOR){X’,(X’ (Y’))∈GEV ⇒ Y’∈GEV}.

Prueba: Supóngase que X’,(X’ (Y’))∈GEV. Por lo que λ ≫ X’,(X’ (Y’)) $\overset{\mathbf{DC}}{\Rightarrow}$ X’ ( (Y’)) $\overset{\mathbf{DCC}}{\Longrightarrow}$ X’ Y’$\overset{\mathbf{B}}{\Rightarrow}$ Y’. Lo cual significa que, λ ≫ Y’ es decir, Y’∈GEV. Por lo tanto, (MP)’ preserva validez. □

Proposición 16. (*Los teoremas de LD son válidos en Gamma-LD*).

(∀Z∈FOR)(Z∈TEO ⇒ Z’∈GEV)

Prueba. Supóngase que Z∈TEO, por lo que existe una demostración de Z. La prueba se realiza por inducción sobre la longitud *n* de la demostración de Z.

**PB**. Paso base. n=0. Esto significa que Z es un axioma, por las proposiciones 13 y 14 se concluye que Z’∈GEV.

**PI**. Paso inductivo. Hipótesis inductiva **HI**: Si la longitud de la demostración de X∈FOR es t con

0 < t < n entonces X’∈GEV. Si la demostración de Z tiene longitud **n** entonces Z, que es el último elemento de la sucesión, o bien es un axioma o se sigue de 2 anteriores utilizando la regla MP. En el primer caso se procede como en el paso base. En el segundo caso se tienen pruebas de X y de X⊃Z, ambas pruebas de longitudes t y k, con t<n y k<n. Aplicando la **HI**, se afirma que X’,X’⊃Z’∈GEV y utilizando la proposición 15, se concluye que Z’∈GEV.

Por el principio de inducción matemática, se concluye que (∀Z∈FOR)(Z∈TEO ⇒ Z’∈GEV). □

Proposición 17. (*Los teoremas de LD son exactamente los gráficos válidos en Gamma-LD*).

(∀Z∈FOR)(Z∈TEO ⇔ Z’∈GEV)

Prueba. Consecuencia directa de las proposiciones 12 y 16.

# Lógica Clásica y Gráficos Existenciales Alfa-LC

*Restringiendo* las reglas de transformación de Gamma-LD a Alfa-LC se obtiene el siguiente conjunto **RTRAC** de reglas primitivas de transformación clásicas, con las cuales se construye el conjunto de ***gráficos existenciales clásicos válidos*** **GCV.** Estas reglas, en esencia, coinciden con las presentadas por Charles Sanders Peirce en [9].

1. **R**$\boldsymbol{\lambda}$: Regla lambda. 2. **B:** Borrado. 3. **E:** Escritura. 4. **DCC**: Doble corte clásico.
5. **DCC**$\lambda$: Doble corte clásico $\lambda$. 6. **I**: Iteración. 7. **D**: Desiteración.
8. **IC**: Iteración clásica. 9. **DC**: Desiteración clásica.

Notación. La ***proposición N**** (*en LI o Alfa-LI*)**,** es el resultado de transformar la ***proposición N*** (*en LD o Gamma-LD*), de tal manera que sólo figuren las restricciones a la *Lógica Clásica* LC.

Simplemente observando las proposiciones ya probadas, y restringiendo en ellas las fórmulas, las reglas de transformación y los axiomas al caso clásico, se concluye que valen las siguientes proposiciones: 2*, 4*, 5*, 7* a 13* y 15* a 17* (a la *proposición N en LD*, le corresponde la *proposición N* en LC*).

Proposición 2* Para X,Y,Z∈**FC**. X⊃Y∈TEO ⇒ ~Y⊃~X∈TEO y ZX⊃ZY∈TEO.

Proposición 4*. $(\forall n\in Z^{+}\cap\{0\})$[X,Y∈**FC** ⇒ $X\bullet G_n(Y) \equiv X\bullet G_n(X\bullet Y)$ y $X\bullet Y^{Rn} \equiv X\bullet(X\bullet Y)^{Rn}$].

Proposición 5*. a. $(\forall n\in Z^{+}\cap\{0\})$[X,Y∈**FC** y X⊃Y∈TEO ⇒ $F_{2n}(X)\supset F_{2n}(Y)$∈TEO ($X^{R2n} \supset Y^{R2n}$∈TEO)].

b. $(\forall n\in Z^{+}\cap\{0\})$[X,Y∈**FC** y X⊃Y∈TEO ⇒ $F_{2n+1}(Y)\supset F_{2n+1}(X)$∈TEO ($Y^{R2n+1} \supset X^{R2n+1}$∈TEO)].

Proposición 7*. Para X,Y∈**Alfa-LC**. a. (∀R∈RTRA)[X∈RP y $X\xRightarrow{R}Y$](∃R'∈RTRA)[Y∈RI y $Y\xRightarrow{R'}X$].

b. (∀R∈RTRA){[X∈$R_k$P y $X\xRightarrow{R}Y$ ⇒ Y∈$R_k$P] y [X∈$RI_k$ y $X\xRightarrow{R}Y$ ⇒ Y∈$R_k$I]}.

Proposición 8*. Para $X_0,X_n$∈**Alfa-LC**. $X_0$>>$X_n$ ⇒ {($X_0$∈RP ⇒ $X_0$>>$X_n$) y ($X_n$∈RI ⇒ $X_n$>>$X_0$)}.

Proposición 9*. TDG. Para X,Y∈**Alfa-LC**. X>>Y ⇒ (X (Y)).

Proposición 10*. TDIG. Para X∈**Alfa-LC**. a. X>>( ) ⇒ (X). b. (X)>>( ) ⇒ X.

Proposición 11*. (∀X,Y∈**FC**)($X'\overset{R}{\Rightarrow}Y'$∈RTRA ⇒ X⊃Y∈TEO).

Proposición 12*. (∀Z∈**FC**)(Z'∈**GCV** ⇒ Z∈TEO).

Proposición 13*. (∀X∈**AxLC**)(X'∈**GCV**).

Proposición 15*. (∀X,Y∈**FC**){X',(X' (Y'))∈**GCV** ⇒ Y'∈**GCV**}.

Proposición 16*. (∀Z∈**FC**)(Z∈TEO ⇒ Z'∈**GCV**).

Proposición 17*. $(\forall Z\in \mathbf{FC})(Z\in TEO \Leftrightarrow Z'\in \mathbf{GCV})$. Es decir, los *teoremas de la lógica proposicional clásica LC* son exactamente los *gráficos existenciales Alfa-LC válidos*.

# Lógica Intuicionista y Gráficos Alfa Intuicionistas

El *sistema deductivo* para la *Lógica Proposicional Intuicionista* **LI** consta de los siguientes axiomas (donde $X,Y,Z\in \mathbf{FI}$):

Ax1.1i. $X\rightarrow(Y\rightarrow X)$
Ax1.2i. $(X\rightarrow(Y\rightarrow Z))\rightarrow((X\rightarrow Y)\rightarrow(X\rightarrow Z))$
Ax1.3i. $X\rightarrow(X\vee Y)$
Ax1.4i. $Y\rightarrow(X\vee Y)$
Ax1.5i. $(X\rightarrow Z)\rightarrow((Y\rightarrow Z)\rightarrow((X\vee Y)\rightarrow Z))$
Ax1.6i. $(X\wedge Y)\rightarrow X$
Ax1.7i. $(X\wedge Y)\rightarrow Y$
Ax1.8. $(X\rightarrow Y)\rightarrow((X\rightarrow Z)\rightarrow(X\rightarrow(Y\wedge Z)))$
Ax1.9i. $X\rightarrow(\neg X\rightarrow Y)$
Ax1.10i. $(X\rightarrow\neg Y)\rightarrow(Y\rightarrow\neg X)$
Ax1.11i. $(X\leftrightarrow Y)\rightarrow(X\rightarrow Y)$
Ax1.12i. $(X\leftrightarrow Y)\rightarrow(Y\rightarrow X)$
Ax1.13i. $(X\rightarrow Y)\rightarrow[(Y\rightarrow X)\rightarrow(X\leftrightarrow Y)]$

Como única *regla de inferencia* se tiene el *Modus Ponens* MPi: de X y $X\rightarrow Y$ se infiere Y.

Definición 7 (*Axiomas de LI*). AxLI ={Ax1.1i, …, Ax1.13i}. Esta presentación de la Lógica Intuicionista es una variante, bastante conocida, de la propuesta original presentada por Arend Heyting en [10].

Proposición 18. (*Conjunción intuicionista*). Para $\underline{X},\underline{Y}\in FI$.

a. $+\underline{X} \equiv \underline{X}$ b. $\underline{X}\bullet\underline{Y} \equiv \underline{X}\wedge\underline{Y}$

Prueba: Parte a. Por AxR se tiene $\neg{\sim}\underline{X}\supset{\sim}{\sim}\underline{X}$, es decir, $+\underline{X}\supset{\sim}{\sim}\underline{X}$, por DN se obtiene $+\underline{X}\supset\underline{X}$. Por AxT se tiene la recíproca $\underline{X}\supset+\underline{X}$. Utilizando Ax1.13 se concluye $+\underline{X}\equiv\underline{X}\in TEO$.

Parte b. Por Ax1.6 se tiene $(\underline{X}\bullet\underline{Y})\supset\underline{X}\in TEO$, por Ax+ resulta $+[(\underline{X}\bullet\underline{Y})\supset\underline{X}]\in TEO$, por MP+ se obtiene $+(\underline{X}\bullet\underline{Y})\supset+\underline{X}\in TEO$, es decir $(\underline{X}\wedge\underline{Y})\supset+\underline{X}\in TEO$, por la parte a, se tiene $+\underline{X}\equiv\underline{X}$, y utilizando EQ y SH se infiere $(\underline{X}\wedge\underline{Y})\supset\underline{X}\in TEO$. De manera similar a partir de Ax1.4 se prueba $(\underline{X}\wedge\underline{Y})\supset\underline{Y}\in TEO$. Utilizando Ax1.8 se deduce $(\underline{X}\bullet\underline{Y})\supset(\underline{X}\wedge\underline{Y})\in TEO$. Al tener los dos condicionales, se concluye que, $(\underline{X}\wedge\underline{Y})\equiv(\underline{X}\bullet\underline{Y})\in TEO$. □

Proposición 19. (*Los axiomas de la lógica proposicional intuicionista son teoremas de LD*).

$(\forall W\in \mathbf{FI})[W\in AxLI \Rightarrow W\in TEO]$

Prueba. Observación. Si $R\in FI$ entonces por las definiciones de LI y de los conectivos intuicionistas se sigue que $R=\underline{a}$ con $\underline{a}\in FAA$ o $R=[S]$ con $S\in FI$, es decir, $R\in FA$, por lo tanto, si $\mathbf{R\in FI}$ entonces $\mathbf{R=\underline{R}}$.
Parte 1. Por Ax1.1 se tiene $X\supset(Y\supset X)\in TEO$, por Ax+ se infiere $+[X\supset(Y\supset X)]\in TEO$, por MP+ resulta $+X\supset+(Y\supset X)\in TEO$, por la proposición 18a se tiene $+X\equiv X$, utilizando EQ y SH se deduce, $X\supset+(Y\supset X)\in TEO$, de nuevo por Ax+ se sigue $+[X\supset+(Y\supset X)]\in TEO$, es decir, $X\rightarrow(Y\rightarrow X)=Ax1.1i\in TEO$.

Parte 2. Supóngase que $X\rightarrow(Y\rightarrow Z)$, $X\rightarrow Y$, X. Por MP se infieren $Y\rightarrow Z$, Y, de nuevo por MP se sigue Z. Aplicando 3 veces TDi, se sigue $(X\rightarrow(Y\rightarrow Z))\rightarrow((X\rightarrow Y)\rightarrow(X\rightarrow Z))=Ax1.2i\in TEO$.

Parte 3. Por Ax1.3 se tiene $X\supset(X\cup Y)\in TEO$, por Ax+ resulta $+[X\supset(X\cup Y)]\in TEO$, por MP+ se obtiene $+X\supset+(X\cup Y)\in TEO$, es decir $+X\supset(X\vee Y)\in TEO$, por la proposición 18a se tiene $+X\equiv X$, y usando EQ y SH se infiere $X\supset(X\vee Y)\in TEO$, por Ax+ resulta $+(X\supset(X\vee Y))\in TEO$, es decir $X\to(X\vee Y)=Ax1.3i \in TEO$.

Parte 4. De manera análoga a la parte 3 se prueba $Y\to(X\vee Y) =Ax1.4i\in TEO$.

Parte 5. Supóngase que $X\to Z$, $Y\to Z$, $X\vee Y$. Utilizando AxR y MP se obtienen $X\supset Z$, $Y\supset Z$, $X\cup Y$, por Ax1.5 y MP se sigue Z. Aplicando 3 veces TDi, se concluye $(X\to Z)\to((Y\to Z)\to ((X\vee Y)\to Z)) =Ax1.5i \in TEO$.

Parte 6. Por Ax1.6 se tiene $(X\bullet Y)\supset X\in TEO$, por Ax+ resulta $+[(X\bullet Y)\supset X]\in TEO$, por MP+ se obtiene $+(X\bullet Y)\supset+X\in TEO$, es decir $(X\wedge Y)\supset+X\in TEO$, por la proposición 18a se tiene $+X\equiv X$, y usando EQ y SH se infiere $(X\wedge Y)\supset X\in TEO$, por Ax+ resulta $+((X\wedge Y)\supset X)\in TEO$, es decir $(X\wedge Y)\to X=Ax1.6i\in TEO$.

Parte 7. De manera análoga a la parte 6 se prueba $(X\wedge Y)\to Y=Ax1.7i\in TEO$.

Parte 8. Supóngase que $X\to Y$, $X\to Z$, X. Utilizando AxR y MP se obtienen $X\supset Y$, $X\supset Z$, X, utilizando Ax1.8 y MP se infiere $Y\bullet Z$, por la proposición 18b se llega a $Y\wedge Z$. Aplicando 3 veces TDi, se concluye $(X\to Y)\to((X\to Z)\to (X\to(Y\wedge Z))) =Ax1.8i\in TEO$.

Parte 9. De manera análoga a la parte 1 se prueba $X\to(\neg X\to Y) =Ax1.9i \in TEO$.

Parte 10. Supóngase que $X\supset\neg Y$, por AxR se tiene que $\neg Y\supset\sim Y$, por SH se infiere $X\supset\sim Y$, por Tras resulta $Y\supset\sim X$. Aplicando TD se obtiene $(X\supset\neg Y)\supset(Y\supset\sim X)\in TEO$. Por proposición 1 se puede afirmar que $+[(X\supset\neg Y)\supset(Y\supset\sim X)]\in TEO$, aplicando MP+ y MP se deduce $+(X\supset\neg Y)\supset+(Y\supset\sim X)\in TEO$, por MP+ se tiene $+(Y\supset\neg X)\supset(+Y\supset+\sim X)\in TEO$, es decir, $+(Y\supset\neg X)\supset(+Y\supset\neg X)\in TEO$, por proposición 1 se infiere $+[+(Y\supset\neg X)\supset(+Y\supset\neg X)]\in TEO$, y por MP+ se sigue $++(Y\supset\neg X)\supset+(+Y\supset\neg X)]\in TEO$, por la proposición 18a se tiene $+Y\equiv Y$ y $++(Y\supset\neg X)\equiv+(Y\supset\neg X)$, entones se deduce que $+(Y\supset\neg X)\supset +(Y\supset\neg X)\in TEO$, lo cual significa que $(Y\to\neg X)\supset(Y\to\neg X)\in TEO$, finalmente por Ax+ se concluye $+[(Y\to\neg X)\supset(Y\to\neg X)]\in TEO$, es decir, $(Y\to\neg X)\to(Y\to\neg X) =Ax1.10i\in TEO$.

Parte 11. Supóngase que $X\leftrightarrow Y$, X. Utilizando AxR y MP se obtiene $X\equiv Y$, aplicando Ax1.11 y MP se infiere $X\supset Y$, y como se tiene X entonces también se tiene Y. Aplicando 2 veces TDi, se concluye $(X\leftrightarrow Y)\to(X\to Y) =Ax1.11i\in TEO$.

Parte 12. De manera análoga a la parte 11 se prueba $(X\leftrightarrow Y)\to(Y\to X) =Ax1.12i\in TEO$.

Parte 13. Supóngase que $X\to Y$, $Y\to X$. Utilizando AxR y MP se obtienen $X\supset Y$, $Y\supset X$, aplicando Ax1.13 y MP se infiere $X\equiv Y$. Si se supone X entonces por Ax1.11 y MP también se tiene Y. Aplicando TDGF (proposición 9b) se obtiene $X\to Y$. De manera similar se obtiene $Y\to X$. Estos dos resultados por EQi significan $X\leftrightarrow Y$. Aplicando 2 veces TDi, se concluye $(X\to Y)\to[(Y\to X)\to(X\leftrightarrow Y)] =Ax1.2i\in TEO$. □

Proposición 20. (*Modus ponens intuicionista vale en LD*). Para $X,Z\in FI$.

MPi. $X\to Z$ y $X \Rightarrow Z$.

Prueba: Supóngase que se tienen $X\to Z$ y X, por lo que resulta $+(X\supset Z)$, por AxR y MP se infiere $X\supset Z$, y como se tiene el antecedente por MP se concluye el consecuente Z. □

Definición 8. (*Teorema de LI*). Para X∈FI. Se dice que X es un *teorema de LI* (X**∈TEOI**) si y solamente si existe una *demostración de X*, es decir, X es la última fórmula de una sucesión finita de fórmulas intuicionistas, tales que cada una de ellas es un axioma o se infiere de dos fórmulas anteriores utilizando la regla de inferencia MPi.

Proposición 21. (*LD incluye LI*).

(∀X∈FI)[X∈TEOI ⇒ X∈TEO]

Prueba. Por inducción matemática sobre la longitud ***n*** de la demostración de X en LI.

**PB**. Paso Base. n=1. Significa que X∈AxLI, por la proposición 19 se infiere que X∈TEO.

**PI**. Paso inductivo. Hipótesis Inductiva **HI**: W∈TEOI ⇒ W∈TEO, donde ***k*** es la longitud de la demostración de W y *0<k<n*. Supóngase que la demostración de X tiene longitud ***n***, se tiene entonces que X es un axioma o X es consecuencia de pasos anteriores utilizando la regla de inferencia MPi. En el primer caso se procede como en el paso base. En el segundo caso se tienen, para alguna fórmula Z, demostraciones de Z y de Z→X, por lo que Z→X∈TEOI y Z∈TEOI, ambas demostraciones de longitud menor que *n*, por **HI** se deduce que Z→X∈TEO y Z∈TEO, aplicando MPi se concluye que Z∈TEO.

Por el principio de inducción matemática, se ha probado que (∀X∈FI)[X∈TEOI ⇒ X∈TEO]. □

# Reglas de transformación y Validez en Alfa-LI

El conjunto **Alfa-LI** de **gráficos existenciales intuicionistas**, ya se definió de la siguiente manera:

1. a∈FAA ⇒ a∈Alfa-LI.
2. λ∈Alfa-LI (λ es el gráfico vacío).
3. X∈Alfa-LI ⇒ **[**X**]**∈Alfa-LI.
4. X,Y∈Alfa-LI ⇒ XY, **[**X **(**Y**)]**, **[(**X**)** **(**Y**)]**∈Alfa-LI.
5. Sólo 1 a 4 determinan Alfa-LI.

Las *reglas de transformación* en Alfa-LI se eligen entre las reglas de Gamma-LD. Estas reglas, en esencia, coinciden con las presentadas por Arnold Oostra en [3].

El conjunto **RTRA-LI** de *reglas primitivas de transformación en Alfa-LI* consta de (donde X,Y,Z∈Alfa-LI):

1. **Rλ**: Regla **λ**.
2. **B:** Borrado.
3. **BL:** Borrado de Lazo. **[**X**]**∈RI ⇒ **[**X **(**Y**)]** $\overset{\mathbf{BL}}{\Rightarrow}$ **[**X**]**
4. **E:** Escritura.
5. **EL:** Escritura de Lazo. **[**X**]**∈RP ⇒ **[**X**]** $\overset{\mathbf{EL}}{\Rightarrow}$ **[**X **(**Y**)]**
6. **R**: Rizo. **[(**X**)]** $\overset{\mathbf{R}}{\Leftrightarrow}$ X
7. **CCR**: Cambio de corte en rizo.

**[**X **(**Y**)]**∈RP ⇒ **[**X **(**Y**)]** $\overset{\mathbf{CCR}}{\Longrightarrow}$ **[**X **[**Y**]]**

**[**X **[**Y**]]**∈RI ⇒ **[**X **[**Y**]]** $\overset{\text{CCR}}{\Longrightarrow}$ **[**X **(**Y**)]**

8. **I**: Iteración.

9. **IdL:** Iteración del Lazo. **[**X **(**Y**)]** $\overset{\text{IdL}}{\Rightarrow}$ **[**X **(**Y**) (**Y**)]**

10. **IF**: Iteración fuerte.

11. **IeL**: Iteración en el lazo. **[**X Y **(**Z**)]** $\overset{\text{IeL}}{\Longrightarrow}$ **[**X Y **(**X Z**)]**

12. **D**: Desiteración.

13. **DdL:** Desiteración del Lazo. **[**X **(**Y**) (**Y**)]** $\overset{\text{DdL}}{\Longrightarrow}$ **[**X **(**Y**)]**

14. **DF**: Desiteración fuerte.

15. **DeL**: Desiteración en el lazo. **[**X Y **(**Z X**)]** $\overset{\text{DeL}}{\Longrightarrow}$ **[**X Y **(**Z**)]**

16. Sólo 1 a 15 determinan RTRA-LI.

Observación. (*seudo-fórmula*). **(**X**)**∉Alfa-LI, **~**X∉FI. Por lo que, no es aplicable la regla de transformación DCC: X ⇒ **((**X**))**, y no es aplicable la regla de inferencia DN: X ⇒ **~~**X. Cuando en Alfa-LI aparezca la *seudo-fórmula* **(**X**)**, sólo aplican las reglas de RTRA-LI, y las reglas derivadas de estas.

Proposición 22. (*Reglas de transformación derivadas*). Para X,Y,Z∈Alfa-LI.

a. **R**λ: Rizo λ. **[( )]**∈GIV

b. **DCI**: Doble corte intuicionista. X∈RP ⇒ X $\overset{\text{DCI}}{\Longrightarrow}$ **[[**X**]]**.

X∈RI ⇒ **[[**X**]]** $\overset{\text{DCI}}{\Longrightarrow}$ X.

c. **TCI**: Triple corte intuicionista. **[[[**X**]]]** $\overset{\text{TCI}}{\Longleftrightarrow}$ **[**X**]**.

d. **IFeL**: Iteración fuerte en el lazo. X[…[Y (Z)]…] $\overset{\text{IFeL}}{\Longrightarrow}$ X[…[Y **(**X Z**)**]…].

e. **DFeL**: Desiteración fuerte en el lazo. X**[**…**[**Y **(**X Z**)]**…**]** $\overset{\text{DFeL}}{\Longrightarrow}$ X[…[Y **(**Z**)**…].

f. **RaN**: Rizo a negación. **[**X **(**[ ]**)]** $\overset{\text{RaN}}{\Longleftrightarrow}$ **[**X**]**

g. **RaD**: Rizo a disyunción.

**[**X **(**Y**) (**Z**)]** $\overset{\text{RaD}}{\Longleftrightarrow}$ **[**X **( [(**Y**) (**Z**)] )]**.

Prueba. Parte a: Por Rλ se tiene λ∈GIV, aplicando R resulta [(λ)]∈GIV, es decir, [( )]∈GEV.

Parte b: X∈RP ⇒ X $\overset{\text{R}}{\Rightarrow}$ [(X)] $\overset{\text{CCR}}{\Longrightarrow}$ [[X]]. X∈RI ⇒ [[X]] $\overset{\text{CCR}}{\Longrightarrow}$ [(X)] $\overset{\text{R}}{\Rightarrow}$ X.

Parte c: [X]∈RP ⇒ [X] $\overset{\text{DCI}}{\Longrightarrow}$ [[[X]]], y [X]∈RP ⇒ [[[X]]] $\overset{\text{CCR}}{\Longrightarrow}$ [[(X)]] $\overset{\text{R}}{\Leftrightarrow}$ [X]. Por lo que, [X]∈RP ⇒ [X] $\overset{\dots}{\Longleftrightarrow}$ [[[X]]]. Además, [X]∈RI ⇒ [[[X]]] $\overset{\text{DCI}}{\Longrightarrow}$ [X], y [X]∈RI ⇒ [X] $\overset{\text{R}}{\Leftrightarrow}$ [([X])] $\overset{\text{CCR}}{\Longrightarrow}$ [[[X]]]. Por lo que, [X]∈RI ⇒ [[[X]]] $\overset{\dots}{\Longleftrightarrow}$ [X]. Por lo tanto, [X] $\overset{\text{TCI}}{\Longleftrightarrow}$ [[[X]]].

Parte d: X[…[…[Y**(**Z**)**]…]…] $\overset{\text{IF}}{\Longrightarrow}$ X[…[… [X Y **(**Z**)**]…]…] $\overset{\text{IeL}}{\Longrightarrow}$ X[…[… [X Y **(**X Z**)**]…]…] $\overset{\text{DF}}{\Longrightarrow}$ X[…[… [Y **(**X Z**)**]…]…].

Parte e: X[…[… [Y **(**X Z**)**]…]…] $\overset{\text{DeL}}{\Longrightarrow}$ X[…[… [X Y **(**Z**)**]…]…] $\overset{\text{DF}}{\Longrightarrow}$ X[…[… [Y **(**Z**)**]…]…]

Parte f: $[X]\in \mathrm{RP} \Rightarrow [X\ ([\ ])] \xLeftrightarrow{\mathbf{Def-\lambda}} [X\ ([\lambda])] \xRightarrow{\mathbf{CCR}} [X\ [[\lambda]]] \xRightarrow{\mathbf{DCI}} [X\ \lambda] \xLeftrightarrow{\mathbf{Def-\lambda}} [X]$, además, $[X]\in \mathrm{RP} \Rightarrow [X] \xRightarrow{\mathbf{EL}} [X\ ([\ ])]$, por lo que, $[X]\in \mathrm{RP} \Rightarrow [X\ ([\ ])] \xLeftrightarrow{\dots} [X]$. Por otro lado, $[X]\in \mathrm{RI} \Rightarrow [X] \xLeftrightarrow{\mathbf{Def-\lambda}} [X\ \lambda] \xRightarrow{\mathbf{DCI}} [X\ [[\lambda]]] \xRightarrow{\mathbf{CCR}} [X\ ([\lambda])] \xLeftrightarrow{\mathbf{Def-\lambda}} [X\ ([\ ])]$, además, $[X]\in \mathrm{RI} \Rightarrow [X\ ([\ ])] \xRightarrow{\mathbf{BL}} [X]$, por lo que, $[X]\in \mathrm{RI} \Rightarrow [X\ ([\ ])] \xLeftrightarrow{\dots} [X]$. Se concluye que, $[X\ ([\ ])] \xLeftrightarrow{\mathbf{RaN}} [X]$.

Parte g: la prueba se presenta después de la proposición 10**. □

En lo que sigue se utilizaran los siguientes resultados de **LI** presentados en [11]. Sean X,Y,Z∈FI (X⇒Y significa que de X se infiere Y. X⇔Y significa que X⇒Y y Y⇒X).

- I∧. *Introducción*. X y Y ⇒ X∧Y.
- E∧. *Eliminación*. X∧Y ⇒ X.
- EQi. *Equivalencia*. X↔Y ⇔ (X→Y)∧(Y→X).
- Trasi.*Transposición*.X→Y⇒¬Y→¬X.
- Impi. *Implicación*. ¬X∨Y ⇒ X→Y ⇒ ¬(X∧¬Y)
- Idem∨. *Idempotencia*. X∨X ⇔ X.
- TDi. *Teorema de deducción* : Si X⇒Y entonces X→Y.
- Idem∧. *Idempotencia*. X∧X ⇔ X.
- DIi. *Demostración indirecta*. X→(Y∧¬Y) ⇒ ¬X.
- Com∧. *Conmutatividad*. X∧Y ⇔ Y∧X.
- Aso∧. *Asociatividad*. X∧(Y∧Z) ⇔ (X∧Y)∧Z.

Notación. La ***proposición N**** (*en LI o Alfa-LI*)**,** es el resultado de transformar la ***proposición N*** (*en LD o Gamma-LD*), de tal manera que sólo figuren las restricciones a la *Lógica Intuicionista* LI.

Simplemente observando las proposiciones ya probadas, restringiendo las fórmulas, las reglas de transformación y los axiomas al caso intuicionista, se concluye que ***siguen siendo válidas*** las siguientes proposiciones: 3**, 5**, 7** a 10**, 12**, 16**, 17**. Sólo deben probarse las proposiciones 2**, 21g, 11**, 14** y 15**. Las proposiciones: 4, 6 y13, no se requerirán en lo que resta de esta sección.

Proposición 2**. (*Trasposición y ampliación)*. Para X,Y,Z∈**FI**. Si X→Y∈TEOI entonces

a. ¬Y→¬X∈TEOI. b. Z∧X→Z∧Y∈TEOI c. (X∧Y→Z)↔(X∧Y→X∧Z)∈TEOI

Prueba: Supóngase que X→Y∈TEOI. Parte a. Resulta de aplicar Trasi en X→Y. Parte b. Supóngase Z∧X, por E∧ se infieren Z y X, por MPi resulta Y, por I∧ se afirma Z∧Y, finalmente por TDi se concluye Z∧X→Z∧Y∈TEOI. Parte c. Supóngase que X∧Y→Z. Supóngase que X∧Y, por E∧ resultan X y Y, por MPi se sigue Z, por I∧ resulta X∧Z, aplicando TDi se obtiene X∧Y→X∧Z. Aplicando TDi de nuevo se concluye que (X∧Y→Z)→(X∧Y→X∧Z)∈TEOI. De manera similar se prueba la recíproca, y por lo tanto, aplicando EQi se deriva (X∧Y→Z)↔(X∧Y→X∧Z)∈TEOI. □

Se restringe la función F de FI en FI de la siguiente manera ($Y, Y_0, Y_1, \dots, Y_{n+1} \in$ FI):

$F_0(Y) = Y_0 \wedge Y$. $F_1(Y) = Y_1 \wedge \neg F_0(Y) = Y_1 \wedge \neg (Y_0 \wedge Y)$. $F_{n+1}(Y) = Y_{n+1} \wedge \neg F_n(Y)$.

Proposición 3**. PIE-ID. (*Principio de inducción estructural para iteración y desiteración*).

$(\forall n\in Z^{+}\cap\{0\})[X,Y\in \mathbf{FI} \Rightarrow X\wedge F_{n}(Y)\leftrightarrow X\wedge F_{n}(X\wedge Y)\in TEO\ (X\wedge Y^{Rn}\leftrightarrow X\wedge (X\wedge Y)^{Rn}\in TEO)]$

Proposición 5**. PIE-RPI. (*Principio de inducción estructural para regiones pares e impares*).

a. $(\forall n\in Z^{+}\cap\{0\})[X,Y\in \mathbf{FI}$ y $X\rightarrow Y\in TEO \Rightarrow F_{2n}(X)\rightarrow F_{2n}(Y)\in TEO\ (X^{R2n}\rightarrow Y^{R2n}\in TEO)]$.

b. $(\forall n\in Z^{+}\cap\{0\})[X,Y\in \mathbf{FI}$ y $X\rightarrow Y\in TEO \Rightarrow F_{2n+1}(Y)\rightarrow F_{2n+1}(X)\in TEO\ (Y^{R2n+1}\rightarrow X^{R2n+1}\in TEO)]$.

Proposición 7**. (*Reversión de las reglas de transformación*). Para $X,Y\in$ **Alfa-LI**.

a. $(\forall R\in$ **RTRA-LI**$)[X\in RP$ y $X\overset{\mathbf{R}}{\Longrightarrow}Y](\exists R'\in RTRA)[Y\in RI$ y $Y\overset{R'}{\Longrightarrow}X]$

b. $(\forall R\in$ **RTRA-LI**$)\{[X\in RP_{k}$ y $X\overset{\mathbf{R}}{\Longrightarrow}Y \Rightarrow Y\in R_{k}P]$ y $[X\in R_{k}I$ y $X\overset{\mathbf{R}}{\Longrightarrow}Y \Rightarrow Y\in R_{k}I]\}$

Proposición 8**. (*Principio de contraposición*). Para $X_{0},X_{n}\in$ **Alfa-LI**.

$X_{0}>>X_{n} \Rightarrow \{(X_{0}\in RP \Rightarrow X_{0}>>X_{n})$ y $(X_{n}\in RI \Rightarrow X_{n}>>X_{0})\}$

Proposición 9**. TDGi. (*Teorema de deducción en Alfa-LI*). Para $X,Y\in$ **Alfa-LI**.

X>>Y ⇒ **[**X **(**Y**)]**

Proposición 10**. TDIGi. (*Teorema de demostración indirecta en Alfa-LI*). Para $X\in$ **Alfa-LI**.

a. X>>**[ ]** ⇒ **[**X**]** b. **[**X]>>**[ ]** ⇒ **[[**X**]]** c. **[[**X**]]**>>**[ ]** ⇒ **[**X]

Prueba de parte g de la proposición 22: [X (Y) (Z)]  X  [( [[ ([ (Y) (Z) ]) ]] )] $\overset{\mathbf{DF}}{\Longrightarrow}$ [(Y) (Z)]  X  [([[ ([ (Y) (Z) ]) ]] )] $\overset{\mathbf{DeL}}{\Longrightarrow}$ [(Y) (Z)]  X [([[ ( ) ]])] $\overset{\mathbf{B}}{\Rightarrow}$ [([[ ( ) ]])] $\overset{\mathbf{R}}{\Leftrightarrow}$ [ [( )] ] $\overset{\mathbf{R}}{\Leftrightarrow}$ [ ]. Se ha probado que [X (Y) (Z)]  X  [( [[ ([ (Y) (Z) ]) ]] )] >> [ ]. Por TDIG, proposición 10**, resulta que [X (Y) (Z)]  >> [X [( [[ ([ (Y) (Z) ]) ]] )]] $\overset{\mathbf{R}}{\Leftrightarrow}$ [X [[ ([ (Y) (Z) ]) ]]] $\overset{\mathbf{DCI}}{\Longrightarrow}$ [X ([(Y) (Z)])]. Por lo tanto, [X (Y) (Z) ]  $\overset{\cdots}{\Rightarrow}$ [X ( [(Y) (Z)] )].

[X ([(Y) (Z)])] X [([[(Y)]])] [([[(Z)]])] $\overset{\mathbf{DF}}{\Longrightarrow}$ [ ([(Y) (Z)])] X [([[(Y)]])] [([[(Z)]])] $\overset{\mathbf{R}}{\Leftrightarrow}$ [(Y) (Z)] X [([[(Y)]])] [([[(Z)]])] $\overset{\mathbf{R}}{\Leftrightarrow}$ [(Y) (Z)] X  [[(Y)]] [[(Z)]] $\overset{\mathbf{R}}{\Leftrightarrow}$ [(Y) (Z)] X  [Y] [Z] $\overset{\mathbf{CCR}}{\Longrightarrow}$ [[Y] [Z]] X  [Y] [Z] $\overset{\mathbf{DF}}{\Longrightarrow}$ [  ] X  [Y] [Z] $\overset{\mathbf{B}}{\Longrightarrow}$ [  ]. Se ha probado que [X ([(Y) (Z)])] X [([[(Y)]])] [([[(Z)]])] >> [  ]. Por TDIGi, proposición 10**, resulta que [X ([(Y) (Z)])] >> [X [([[(Y)]])] [([[(Z)]])]] $\overset{\mathbf{R}}{\Leftrightarrow}$ [X [[(Y)]] [[(Z)]]] $\overset{\mathbf{DCI}}{\Longrightarrow}$ [X (Y) (Z)]. Por lo tanto, [X ( [ (Y) (Z) ] )] $\overset{\cdots}{\Rightarrow}$ [X (Y) (Z) ]. Se concluye finalmente que, [X ( [ (Y) (Z) ] )] $\overset{\mathbf{RaD}}{\Longleftrightarrow}$ [X (Y) (Z) ]. □

Proposición 11**. (*Validez en LI de las reglas primitivas de transformación*). Para $X,Y\in$ **FI**.

$(\forall R\in RTRA\text{-}LI)(\exists T\in TEOI)[X'\overset{\mathbf{R}}{\Rightarrow}Y' \Rightarrow T=X\supset Y]$

Prueba. Regla 1. Rλ: Regla lambda. $\lambda\in GIV$.

$\lambda$=Ax1.1i' y Ax1.1i$\in$TEOI.

Regla 2. B: Borrado. X'Y'$\in$RP $\Rightarrow$ (X'Y' $\overset{\mathbf{B}}{\Rightarrow}$ X') y (X'Y' $\overset{\mathbf{B}}{\Rightarrow}$ Y').

Regla 4. E: Escritura. X'Y'$\in$RI $\Rightarrow$ (X' $\overset{\mathbf{E}}{\Rightarrow}$ X' Y') y (Y' $\overset{\mathbf{E}}{\Rightarrow}$ X'Y').

Por Ax1.6i y Ax1.7i se tienen $(X\wedge Y)\rightarrow X\in$TEOI y $(X\wedge Y)\rightarrow Y\in$TEOI, por la proposición 5** se concluyen que $(X\wedge Y)^{R2n}\rightarrow Y^{R2n}\in$TEOI y $(X\wedge Y)^{R2n}\rightarrow X^{R2n}\in$TEOI, $Y^{R2n+1}\rightarrow(X\wedge Y)^{R2n+1}\in$TEOI y $X^{R2n+1}\rightarrow (X\wedge Y)^{R2n+1}\in$TEOI.

Regla 3 BL: borrado de lazo. X'$\in$RI $\Rightarrow$ [X'] $\overset{\mathbf{EL}}{\Longrightarrow}$ [X' (Y')].

Regla 5 EL: escritura de lazo. X'$\in$RP $\Rightarrow$ [X' (Y')] $\overset{\mathbf{BL}}{\Longrightarrow}$ [X'].

Por Ax1.9i se tiene $\neg X\rightarrow(X\rightarrow Y)\in$TEOI, aplicando la proposición 5** se concluye que $(\neg X)^{R2n}\rightarrow (X\rightarrow Y)^{R2n}\in$TEOI y $(X\rightarrow Y)^{R2n+1}\rightarrow(\neg X)^{R2n+1}\in$TEOI.

Regla 6. R: Rizo. [$\lambda$ (X')] $\overset{\mathbf{R}}{\Leftrightarrow}$ X'.

Por Ax1.1i se tiene $X\rightarrow(((Y\wedge Y)\rightarrow Y)\rightarrow X)\in$TEOI. Supóngase que $((Y\wedge Y)\rightarrow Y)\rightarrow X$, por Ax1.6i se tiene $(Y\wedge Y)\rightarrow Y)$, por MPi se infiere X, por TDi se obtiene que $(((Y\wedge Y)\rightarrow Y)\rightarrow X)\rightarrow X\in$TEOI y como se tiene la recíproca entonces $(((Y\wedge Y)\rightarrow Y)\rightarrow X)\leftrightarrow X\in$TEOI. Aplicando la proposición 5** se concluye que $(X)^{R2n}\leftrightarrow ((Y\wedge Y)\rightarrow Y)\rightarrow X)^{R2n}\in$TEOI y $((Y\wedge Y)\rightarrow Y)\rightarrow X)^{R2n+1}\leftrightarrow X^{R2n+1}\in$TEOI, es decir, $((Y\wedge Y)\rightarrow Y)\rightarrow X)^{Rn}\leftrightarrow X^{Rn}\in$TEOI.

Regla 7. CCR: Cambio de corte en rizo. [X' (Y')]$\in$RP $\Rightarrow$ [X' (Y')] $\overset{\mathbf{CCR}}{\Longrightarrow}$ [X' [Y']] y [X' (Y')]$\in$RI $\Rightarrow$ [X' [Y']] $\overset{\mathbf{CCR}}{\Longrightarrow}$ [X' (Y')].

Por Impi se tiene $(X\rightarrow Y)\rightarrow\neg(X\wedge\neg Y)\in$TEOI, aplicando la proposición 5** se concluye que $(X\rightarrow Y)^{R2n}\rightarrow (\neg(X\wedge\neg Y))^{R2n}\in$TEOI y $(\neg(X\wedge\neg Y))^{R2n+1}\rightarrow(X\rightarrow Y)^{R2n+1}\in$TEOI.

Regla 8. **I**: Iteración. X'$\overset{\mathbf{I}}{\Rightarrow}$X'X'.

Regla 12. **D**: Desiteración. X'X'$\overset{\mathbf{D}}{\Rightarrow}$X'.

Por Idem$\wedge$ se tiene $(X\wedge X)\leftrightarrow X\in$TEOI aplicando la proposición 5** se concluye que $(XX)^{R2n}\leftrightarrow X^{R2n}\in$TEOI y $(XX)^{R2n+1}\leftrightarrow X^{R2n+1}\in$TEOI, es decir, $(XX)^{Rn}\leftrightarrow X^{Rn}\in$TEOI.

Regla 10. **IF**: Iteración fuerte. X'[…[Y']…] $\overset{\mathbf{IF}}{\Longrightarrow}$ X'[…[X'Y']…].

Regla 14. **DF**: Desiteración fuerte. X'[…[X'Y']…] $\overset{\mathbf{DF}}{\Longrightarrow}$ X'[…[Y']…].

Por la proposición 3** se tiene que, para $n\in Z^{+}\cap\{0\}$, $X\wedge Y^{Rn}\leftrightarrow X\wedge(X\wedge Y)^{Rn}$.

Regla 11. IeL: Iteración fuerte en el lazo. [X' Y' **(**Z'**)**] $\overset{\mathbf{IeL}}{\Longrightarrow}$ [X' Y' **(**X' Z'**)**].

Regla 15. DeL: Desiteración en el lazo. [X' Y' **(**X' Z'**)**] $\overset{\mathbf{DeL}}{\Longrightarrow}$ [X' Y' **(**Z'**)**].

Por la proposición 2c** se tiene $(XY\to Z)\leftrightarrow(XY\to XZ)\in$TEOI, por la proposición 5** resultan $(XY\to Z)^{R2n}\leftrightarrow(XY\to XZ)^{R2n}\in$TEOI y $(XY\to XZ)^{R2n+1}\leftrightarrow(XY\to Z)^{R2n+1}\in$TEOI, y así, $(XY\to Z)^{Rn}\leftrightarrow(XY\to XZ)^{Rn}\in$TEOI.

Regla 9. IdL: Iteración de lazo. **[**X' **(**Y'**)]** $\overset{\textbf{IdL}}{\Longrightarrow}$ [X' **(**Y'**) (**Y'**)]**.

Regla 13. DdL: Desiteración de lazo. [X' **(**Y'**) (**Y'**)]** $\overset{\textbf{DdL}}{\Longrightarrow}$ **[**X' **(**Y'**)].**

Por idemp$\vee$ se tiene $Y\leftrightarrow(Y\vee Y)\in$TEOI, por lo que $(X\to Y)\leftrightarrow(X\to(Y\vee Y))\in$TEOI. Aplicando la proposición 5** se concluye que $(X\to Y)^{R2n}\leftrightarrow(X\to(Y\vee Y))^{R2n}\in$TEOI y $(X\to Y)^{R2n+1}\leftrightarrow(X\to(Y\vee Y))^{R2n+1}\in$TEOI, es decir, $(X\to Y)^{Rn}\leftrightarrow(X\to(Y\vee Y))^{Rn}\in$TEOI, lo cual significa que [X' **(**Y'**)]** $\overset{\cdots}{\Longleftrightarrow}$ [X' **( [(**Y'**) (**Y'**)] )]** $\overset{\textbf{RaN}}{\Longleftrightarrow}$ [X' **(**Y'**) (**Y'**)]**$^{\$}$.

$^{\$}$: Observar en la proposición 22g, que la regla RaD puede ser aplicada, puesto que, para ser garantizada, se requieren las reglas primitivas 14, 15, 2, 6 y 7, las cuales ya se garantizaron en los párrafos precedentes de la prueba de esta proposición (11**). □

Proposición 12**. (*Los gráficos válidos en Gamma-LD son teoremas en LD*).

$(\forall Z\in FOR)(Z'\in GIV \Rightarrow Z\in TEOI)$

Prueba: Si $\lambda>>X_n$' entonces existen $R_1$, …, $R_n\in$RTRA-LI, y existen $\lambda$, $X_2$', …, $X_{n-1}$'$\in$Alfa-LI, tales que $\lambda R_1X_1'R_2X_2'\ldots X_{n-1}'R_nX_n'$, y se dice que la *longitud* de la transformación de $\lambda>>X_n$' es **n**.

La prueba se realiza por inducción sobre la longitud *n* de la transformación.

**PB**. Paso base. *n*=1. Significa que solo se aplicó una de las reglas primitivas, y $R_1$ tiene que ser de la forma $\lambda \overset{\textbf{R1}}{\Longrightarrow} X_1$' . Por la proposición 11** se infiere que $Ax1.1\to X_1\in$TEO, luego por MPi $X_1\in$TEO.

**PI**. Paso inductivo. Hipótesis inductiva **HI**: $(\forall n>1)[\lambda>>X_n' \Rightarrow X_n\in TEOI$. Si $\lambda>>X_{n+1}$' entonces $\lambda R_1X_1'R_2X_2'\ldots X_{n-1}'R_nX_n'R_{n+1}X_{n+1}'$, es decir, $\lambda R_1X_1'R_2X_2'\ldots X_{n-1}'R_nX_n'$ y $X_n'R_{n+1}X'_{n+1}$, por lo que, $\lambda>>X_n$' y $X_n'R_{n+1}X_{n+1}'$. Aplicando **HI** y TDG proposición 11** se infieren $X_n\in$TEOI y $X_n\to X_{n+1}\in$TEOI, aplicando MPi se concluye que $X_{n+1}\in$TEOI. Por el principio de inducción matemática se tiene que $(\forall Z\in FOR)$ $(Z'\in GEV \Rightarrow Z\in TEO)$. □

Proposición 14**. (*AxLI válidos en Alfa-LI*).

$(\forall X\in$**LI**$)(X\in TEOI \Rightarrow X'\in GIV)$

Prueba. Ax1.1: $\lambda \overset{\textbf{R}}{\Rightarrow}$ [(λ)] $\overset{\textbf{R}}{\Rightarrow}$ [( [( )] )] $\overset{\textbf{E}}{\Rightarrow}$ [X'( [Y' ( )] )] $\overset{\textbf{IeL}}{\Rightarrow}$ [X'( [Y' (X')] )]. Por lo tanto, (Ax1.1)'$\in$GEV.

Ax1.2: [X'([Y'(Z')])] [X'(Y')] X'$\overset{\textbf{DF}}{\Longrightarrow}$ [X'([Y'(Z')])] [(Y')] X' $\overset{\textbf{R}}{\Rightarrow}$ [X'([Y'(Z')])] Y' X'$\overset{\textbf{DF}}{\Longrightarrow}$ [X'([(Z')])] Y' X'$\overset{\textbf{R}}{\Rightarrow}$ [X'(Z')] Y' X' $\overset{\textbf{DF}}{\Longrightarrow}$ [(Z')] Y' X'$\overset{\textbf{R}}{\Rightarrow}$ Z' Y' X'$\overset{\textbf{B}}{\Rightarrow}$ Z'. Aplicando TDGi (proposición 9**) 3 veces resulta que (Ax1.2)'$\in$GIV.

Ax1.3 y Ax1.4: $\lambda$ $\overset{\mathbf{R}}{\Rightarrow}$ [($\lambda$)] $\overset{\mathbf{R}}{\Rightarrow}$ [( [( )] )] $\overset{\mathbf{E}}{\Rightarrow}$ [X’( [(Y’) ( )] )] $\overset{\mathbf{IeL}}{\Longrightarrow}$ [X’( X’ [(Y’) ( )] )] $\overset{\mathbf{IF}}{\Rightarrow}$ [X’(X’[ X’ (Y’) ( )])] $\overset{\mathbf{IeL}}{\Longrightarrow}$ [X’( X’ [ X’ (Y’) (X’)] )] $\overset{\mathbf{DF}}{\Longrightarrow}$ [X’( X’ [(Y’) (X’)] )] $\overset{\mathbf{DeL}}{\Longrightarrow}$ [X’( [(Y’) (X’)] )]. Por lo tanto, (Ax1.3)’, (Ax1.4)’$\in$GIV.

Ax1.5: [X’(Z’)] [Y’(Z’)] [(X’)(Y’)] [[(Z’)]] $\overset{\mathbf{R}}{\Rightarrow}$ [X’ (Z’)] [Y’ (Z’)] [(X’)(Y’)] [Z’] $\overset{\mathbf{CCR}}{\Longrightarrow}$ [X’[Z’]] [Y’[Z’]] [[X’] [Y’]] [Z’] $\overset{\mathbf{DF}}{\Rightarrow}$ [X’] [Y’] [[X’] [Y’]] [Z’] $\overset{\mathbf{DF}}{\Rightarrow}$ [X’] [Y’] [ ] [Z’] $\overset{\mathbf{B}}{\Rightarrow}$ [ ]. Por TDIGi (proposición 10**) se infiere que [X’(Z’)] [Y’(Z’)] [(X’)(Y’)] >> [[[(Z’)]]] $\overset{\mathbf{TCI}}{\Longleftrightarrow}$ [(Z’)] $\overset{\mathbf{R}}{\Rightarrow}$ Z’. Por TDGi (proposición 9**) se infiere que [X’(Z’)] [Y’(Z’)] >> [ [(X’)(Y’)] (Z’)]. Por lo tanto, (Ax1.5)’$\in$GIV.

Ax1.6 y Ax1.7: $\lambda$ $\overset{\mathbf{R}}{\Rightarrow}$ [( )] $\overset{\mathbf{E}}{\Rightarrow}$ [X’Y’( )] $\overset{\mathbf{IeL}}{\Longrightarrow}$ [X’Y’(X’)]. Por lo tanto, (Ax1.6)’,(Ax1.7)’$\in$GIV.

Ax1.8: [X’(Y’)] [X’(Z’)] X’ $\overset{\mathbf{DF}}{\Longrightarrow}$ [(Y’)] [(Z’)] X’ $\overset{\mathbf{R}}{\Rightarrow}$ Y’ Z’ X’ $\overset{\mathbf{B}}{\Rightarrow}$ Y’ Z’. Aplicando TDGi 2 veces resulta que (Ax1.8)’$\in$GIV.

Ax1.9: $\lambda$ $\overset{\mathbf{R}}{\Rightarrow}$ [($\lambda$)] $\overset{\mathbf{R}}{\Rightarrow}$ [( [( )] )] $\overset{\mathbf{E}}{\Rightarrow}$ [X’( [( ) (Y’)] )] $\overset{\mathbf{CCR}}{\Longrightarrow}$ [X’( [ [ ] (Y’) ] )] $\overset{\mathbf{IeL}}{\Longrightarrow}$ [X’( X’[ [ ] (Y’) ] )] $\overset{\mathbf{IF}}{\Rightarrow}$ [X’(X’ [[X’] (Y’)])] $\overset{\mathbf{DeL}}{\Longrightarrow}$ [X’ ([[X’] (Y’)])] . Por lo tanto, (Ax1.9)’$\in$GIV.

Ax1.10: [X’ ( [Y’] )] Y’ $\overset{\mathbf{IF}}{\Rightarrow}$ [X’ ( [Y’] ) Y’] Y’ $\overset{\mathbf{IeL}}{\Longrightarrow}$ [X’ ( [Y’] Y’) Y’] Y’ $\overset{\mathbf{DF}}{\Longrightarrow}$ [X’ ( [ ] Y’) Y’] Y’ $\overset{\mathbf{DeL}}{\Longrightarrow}$ [X’ ( [ ] ) Y’] Y’ $\overset{\mathbf{DF}}{\Longrightarrow}$ [X’ ( [ ] )] Y’ $\overset{\mathbf{B}}{\Rightarrow}$ [X’ ( [ ] )] $\overset{\mathbf{RaN}}{\Longrightarrow}$ [X’]. Aplicando TDGi 2 veces resulta (Ax1.10)’$\in$GIV.

Ax1.11: (X$\leftrightarrow$Y)’ $\overset{\mathbf{Definición}}{\Longrightarrow}$ (X$\rightarrow$Y)’ (Y$\rightarrow$X)’ $\overset{\mathbf{B}}{\Rightarrow}$ (X$\rightarrow$Y)’. Por lo tanto, (Ax1.11)’$\in$ GIV.

Ax1.12: (X$\leftrightarrow$Y)’ $\overset{\mathbf{Definición}}{\Longrightarrow}$ (X$\rightarrow$Y)’ (Y$\rightarrow$X)’ $\overset{\mathbf{B}}{\Rightarrow}$ (Y$\rightarrow$X)’. Por lo tanto, (Ax1.12)’$\in$ GIV.

Ax1.13: (X$\rightarrow$Y)’ (Y$\rightarrow$X)’ $\overset{\mathbf{Definición}}{\Longrightarrow}$ (X$\leftrightarrow$Y)’. Aplicando TDGi 2 veces resulta que (Ax1.13)’$\in$GIV. □

Proposición 15**. (*MPi)’ preserva validez en Alfa-LI*).

X’$\in$**GIV** y (X$\rightarrow$Y)’$\in$GIV $\Rightarrow$ Y’$\in$GIV.

Prueba: Supóngase que X’$\in$GIV y (X$\rightarrow$Y)’$\in$GIV. Por lo que $\lambda$ >> X’ [X’ (Y’)] $\overset{\mathbf{DF}}{\Longrightarrow}$ X’ [(Y’)] $\overset{\mathbf{R}}{\Rightarrow}$ X’ Y’ $\overset{\mathbf{B}}{\Rightarrow}$ Y’. Lo cual significa que, $\lambda$ >> Y’, es decir, Y’$\in$GIV. Por lo tanto, (MPi)’ preserva validez. □

Proposición 16**. (*Los teoremas de LI son válidos en Alfa-LI*).

($\forall$Z$\in$**FI**)(Z$\in$TEO $\Rightarrow$ Z’$\in$GIV)

Proposición 17**. (*Los teoremas de LI son exactamente los gráficos válidos en Alfa-LI*).

($\forall$Z$\in$**FI**)(Z$\in$TEO $\Leftrightarrow$ Z’$\in$GIV)

Es decir, los teoremas de la lógica proposicional intuicionista LI son exactamente los gráficos existenciales válidos en Alfa-LI. □

# Consistencia de LD y de Gamma-LD

Sea { }” la función de traducción de Gamma-LD en Alfa-LC, definida así:

1. [X]” = (X)” = (X”) 2. {X Y}” = X” Y” 3. $\lambda$” = $\lambda$

Proposición 23. (*Las reglas de Gamma-LD se transforman en reglas de Gamma-LC*).

X∈RTRA ⇒ X”∈RTRAC.

Prueba.

1. {R$\lambda$}” = R$\lambda$. 2. B” = B. 3. E” = E. 4. DCC” = DCC

5. CC” = {[X] $\overset{\textbf{CC}}{\Rightarrow}$ (X)}’’ = (X’’) $\xRightarrow{\textbf{Identidad}}$ (X’’). 6. DCMGEV” = DCC.

7. DCMF” = DCC. 8. I” = I. 9. D” = D. 10. IC” = IC.

11. DC” = DC. 12. IF” = IC. 13. DF” = DC. □

Proposición 24. (*Los gráficos de Gamma-LD se transforman en gráficos de Gamma-LC*).

Para $X_n$∈ Gama-LD, $X_n$∈GEV ⇒ $X_n$”∈GCV.

Prueba. Si $X_n$∈GEV entonces $\lambda$>>$X_n$, es decir, existen $X_i$∈Gama-LD y $R_i$∈RTRA con 0<i<n+1, tales que $\lambda R_1X_1R_2X_2$ … $X_{n-1}R_nX_n$. Como $R_i$∈RTRA por la proposición 23 resulta que $R_i$”∈RTRAC, además, $X_i$’’∈Alfa-LC obteniéndose $\lambda R_1$”$X_1$”$R_2$”$X_2$” … $X_{n-1}$”$R_n$”$X_n$”, por lo que $\lambda$>>$X_n$”, es decir, $X_n$”∈GCV. □

Proposición 25. (*Consistencia de LD y de Gamma-LD*).

1. ( )∉GEV 2. [ ]∉GEV 3. LD es consistente. 4. LI es consistente.

Prueba. Parte a. Si ( )∈GEV, entonces por la proposición 24 se deduce ( )”∈GCV, por lo que ( )∈GCV, por la proposición 17* esto significa que existe Z∈FC tal que Z•~Z es un teorema de LC, lo cual es imposible, puesto que LC es consistente. Por lo tanto, ( )∉GEV, dicho de otra manera, *Gamma-LD es consistente*.

Parte b. Si [ ]∈GEV entonces por la regla CC se sigue ( )∈GCV, lo cual contradice la parte a. Por lo tanto, [ ]∉GEV.

Parte c. Si LD es inconsistente, entonces existe X∈FOR, tal que X•~X∈TEO, por la proposición 17 se infiere que (X•~X)’∈GEV, es decir, ( )∈GEV, lo cual contradice la parte a. Por lo tanto, se garantiza la *consistencia de LD*.

Parte d. Si LI es inconsistente, entonces existe X∈FOR, tal que X•¬X∈TEOI, por la proposición 21 se infiere que X•¬X∈TEO, es decir, LD es inconsistente, lo cual contradice la parte c. Por lo tanto, se garantiza la *consistencia de LI*. □

# Conclusiones

Conclusión 1. (*Restricción de Gamma-LD a Alfa-LC)*. Por la proposición 17*, si los gráficos existenciales de Gamma-LD se restringen Alfa-LD, resulta que el conjunto de gráficos válidos, en el lenguaje {~, ⊃, •, ∪, ≡}, coincide con los gráficos existenciales *Alfa-Peirce* presentados en [1], es decir, con el conjunto de teoremas del cálculo proposicional clásico. Por lo que, los teoremas del cálculo proposicional clásico son válidos en LD. □

Conclusión 2. (*Restricción de Gamma-LD a Alfa-LI)*. Por la proposición 17**, si los gráficos existenciales de Gamma-LD se restringen Alfa-LI, en el lenguaje {¬, ∧, ∨ →, ↔}, resulta que el conjunto de gráficos válidos coincide con los gráficos existenciales *Alfa intuicionistas* presentados por Arnold Oostra en [3], es decir, con el conjunto de teoremas del cálculo proposicional intuicionista. Por lo que, los teoremas del cálculo proposicional intuicionista son válidos en LD. □

Platón en uno de sus diálogos, *Crátilo* [12], define la verdad como “El discurso, que dice las cosas como son, es verdadero; y el que las dice como no son, es falso”. En el libro IV de la Metafísica (1011b) [4], Aristóteles define el concepto de verdad de la siguiente manera “decir de *lo que es* que es, y de *lo que no es* que no es, es *lo verdadero*; decir de *lo que es* que no es, y de *lo que no es* que es, es *lo falso*”.

Tomando estas definiciones como axiomas, en [13] se construye el sistema *Lógica Básica para la Verdad Aristotélica* **LBVA** (o LBPcPo), resultando que LD es una extensión propia de este. Por lo que, en el sistema LD se pueden modelar estas definiciones interpretando +X como ‘X es verdadero’, ¬X como ‘X es falso’, *X como ‘decir X’, X como ‘X es’ y ~X como ‘X no es’. La definición de *Verdad Aristotélica* estaría codificada por la fórmula (V•*V)≡+V, es decir, por el gráfico **[V[([(V)]) ([V])] ([(V)])] [V (V [([(V)]) ([V])])]**, y la definición de *Falsedad Aristotélica* estaría codificada por la fórmula (~F•*F)≡¬F, es decir, por el gráfico **[(F) [([(F)]) ([F])] ([F])] [[F] ((F) [([(F)]) ([F])])]**.

Conclusión 3. (*Verdad y Falsedad Aristotélica en LD*). Para V,F∈FOR.

a. [(V•*V) ≡ +V]’∈GEV. b. [(~F•*F) ≡ ¬F]’∈GEV.

Prueba: Parte a. Para Z∈FOR, es decir, Z’∈Gamma-LD

{V•*V}’ $\overset{\text{EQ}}{\Longleftrightarrow}$ **(**V•(+V∪¬V)**)**’ $\overset{\text{Traducción}}{\Longleftrightarrow}$ V’ **(([(**V’**)]) ([**V’**]))** $\overset{\text{CC}}{\Longrightarrow}$ V’ **(([(**V’**)]) ((**V’**)))** $\overset{\text{DCC}}{\Longrightarrow}$ V’ **(([(**V’**)])** V’**)** $\overset{\text{D}}{\Longrightarrow}$ V’ **(([(**V’**)]))** $\overset{\text{DCC}}{\Longrightarrow}$ V’ **[(**V’**)]** $\overset{\text{B}}{\Rightarrow}$ **[(**V’**)]** $\overset{\text{Traducción}}{\Longleftrightarrow}$ **{**+V}’. Aplicando TDG proposición 9 resulta que, λ >> **(**{V•*V}’ **({**+V}’**))**, es decir, λ >> {(V•*V)⊃ **(**+V)}’.

{+V}’ $\overset{\text{Traducción}}{\Longleftrightarrow}$ **[(**V’**)]** $\overset{\text{I}}{\Rightarrow}$ **[(**V’**)] [(**V’**)]** $\overset{\text{CC}}{\Longrightarrow}$ **((**V’**)) [(**V’**)]** $\overset{\text{I}}{\Rightarrow}$ **((**V’**)) ((**V’**)) [(**V’**)]** $\overset{\text{DCC}}{\Longrightarrow}$ V’ **((**V’**)) (([(**V’**)]))** $\overset{\text{CC}}{\Longrightarrow}$ V’ **([**V’**]) (([(**V’**)]))** $\overset{\text{I}}{\Rightarrow}$ V’ **([**V’**]) (([**V’**]) ([(**V’**)]))** $\overset{\text{B}}{\Rightarrow}$ V’ **(([**V’**]) ([(**V’**)]))** $\overset{\text{Traducción}}{\Longleftrightarrow}$ **{**V•(+V∪¬V)}’. Aplicando TDG resulta λ >> **(**{+V}’ **({**V•(+V∪¬V)}’**))**, es decir, λ >> {+V⊃(V•*V)}’.

De los dos resultados obtenidos se sigue, λ >> {+V ≡ (V•*V)}’, es decir, {+V≡(V•*V)}’∈GEV, utilizando la proposición 17 se concluye que +V≡(V•*V)∈TEO. De manera similar se prueba la parte b. □

Muchas *paradojas lógicas* involucran los conceptos de verdad o falsedad, por ejemplo, la siguiente variante de la paradoja del mentiroso en [5] y [6]. Considérese la situación en la cual se tiene una oración que dice:

Esta oración es falsa

Cuando se identifican en la lógica clásica, ser el caso con verdadero (Z) y no ser el caso con falso (~Z), entonces se tiene la paradoja: si la oración es el caso (Z) entonces resulta que también es falsa (~Z), y si la oración es falsa (~Z) entonces resulta que es el caso (Z). Es decir, (Z⊃~Z)•(~Z⊃Z). Resultando Z•~Z. Como la lógica clásica no soporta las contradicciones, resulta inútil en este caso.

Cuando se identifican en la lógica intuicionista, ser el caso como verdadero (Z), y no ser el caso con falso (¬Z), entonces se tiene la paradoja: si la oración es el caso (Z) entonces resulta que también es falsa (¬Z), y si la oración es falsa (¬Z) entonces resulta que es el caso (Z). Es decir, (Z→¬Z)∧(¬Z→Z). Resultando Z∧¬Z. Como la lógica intuicionista no soporta las contradicciones, resulta inútil en este caso.

Cuando se identifican en la lógica intuicionista, ser el caso (Z), la doble negación como verdadero (¬¬Z), y no ser el caso con falso (¬Z), entonces se tiene la paradoja: si la oración es verdadera (¬¬Z) entonces resulta que también es falsa (¬Z), y si la oración es falsa (¬Z) entonces resulta que es verdadera (¬¬Z). Es decir, (¬¬Z→¬Z)∧(¬Z→¬¬Z). Resultando ¬Z∧¬¬Z. Como la lógica intuicionista no soporta las contradicciones, resulta inútil en este caso.

Conclusión 4. (*Solución a una versión de la paradoja del mentiroso*).

Con el fin de representar simbólicamente en la Lógica Doble LD, la oración que dice ‘*esta oración es falsa*’, se utiliza la siguiente interpretación:

W significa ‘*lo que dice la oración es el caso*’. ~W significa ‘*lo que dice la oración no es el caso*’.

¬W significa ‘*la oración W es falsa*’. +W significa ‘*la oración W es verdadera*’.

La oración que dice ‘*esta oración es falsa*’, queda entonces representada por la fórmula W≡¬W.

Como consecuencias en LD se tienen los siguientes resultados, donde W∈FAA:

1. W≡¬W no genera contradicciones.

2. A partir de W≡¬W se infieren como conclusiones válidas: ~W, ~¬W, ~+W y ~*W.

Lo anterior significa que, W no es el caso, no es falsa y no es verdadera, lo cual significa que W no está bien fundada, es decir, W no dice nada, W no tiene valor de verdad.

Prueba:

Parte 1. Sea W∈FAA. Supóngase que en LD, se tiene que (W≡¬W)⊃(Z•~Z) para alguna Z∈FOR, es decir, **(**W’**(**[W’]**))(**[W’]**(**W’**))** >> **( )**, lo cual por TDIG (proposición 10a) implica **((**W’**(**[W’]**)) (**[W’]**(**W’**)))**.

Si además se tuviese **(**W’**) ([**W’**])**, por DC en el resultado previo, se infiere **((**W’**([**W’**])) ([**W’**]))**, de nuevo por DC se sigue **((**W’**([**W’**])))**, aplicando DCC resulta W’**([**W’**])**, por regla B se deduce W’, y como se tiene **(**W’**)**, por DC se obtiene **( )**. Por TDIG (proposición 10a) se afirma que **((**W’**) ([**W’**]))**, es decir, **(**W’**)** >> **[**W’**]**.

Por la regla CC se sabe que **[**W’**]** >> **(**W’**)**, y como ya se probó **(**W’**)** >> **[**W’**]**, se concluye que los cortes clásico e intuicionista, se pueden intercambiar en cualquier región par o impar. Por otro lado, la regla DCC dice que **((**W’**))** >> W’, aplicando TDG (proposición 9a) se llega a que **(((**W’**)) (**W’**))**∈GEV, intercambiando cortes se puede asegurar que **[[[**W’**]] (**W’**)]**∈GEV, y como W∈FAA, se deduce **[[[**W’**]] (**W’**)]**∈GIV, se concluye de esta manera que ¬¬W→W∈TEOI, lo cual no es el caso. Por lo tanto, en LD, no se infiere (W≡¬W)⊃(Z•~Z) para alguna Z∈FOR, dicho de otra manera, W≡¬W no genera una contradicción, lo cual significa que, en LD, no se tiene esta versión de la paradoja del mentiroso.

Parte 2. {W≡¬W}’ $\overset{\textbf{EQ}}{\Longleftrightarrow}$ {(W⊃¬W)•(W⊃¬W)}’ $\overset{\textbf{Traducción}}{\Longleftrightarrow}$ **(**W’ **([**W’**])) ((**W’**) [**W’**])** $\overset{\textbf{CC}}{\Longrightarrow}$ **(**W’ **([**W’**])) ([**W’**] [**W’**])** $\overset{\textbf{D}}{\Longrightarrow}$ **(**W’ **([**W’**])) ([**W’**])** $\overset{\textbf{D}}{\Longrightarrow}$ **(**W’**) ([**W’**])** $\overset{\textbf{I}}{\Rightarrow}$ **(**W’**) ([**W’**]) (**W’**)** $\overset{\textbf{DCC}}{\Longrightarrow}$ **(**W’**) ([**W’**]) (((**W’**)))** $\overset{\textbf{CC}}{\Longrightarrow}$ **(**W’**) ([**W’**]) ([(**W’**)])** $\overset{\textbf{Traducción}}{\Longleftrightarrow}$ **{**~W • ~¬W • ~+W}’ $\overset{\textbf{N∪}}{\Longleftrightarrow}$ **{**~W • ~(¬W ∪ +W)}’ $\overset{\textbf{Def*}}{\Longleftrightarrow}$ {~W • ~*W}’.

Conclusión 5. (*Ubicación de LD en la jerarquía de LBVA*). En la jerarquía de lógicas basadas en la *Lógica Básica Paraconsistente y Paracompleta* **LB** (presentada en [14]), se tiene que: la Lógica Doble LD es independiente de la Lógica Básica para la Verdad y la Falsedad **LBVF** (presentada en [15]), y además, LD se encuentra acotada por la *Lógica de las Tautologías* **LT** (presentada en [16]) y la Lógica Básica para la Verdad Aristotélica **LBVA** (presentada en [13]).

Preguntas abiertas. (*Existencia de nuevos sistemas de gráficos existenciales*)

1) ¿Existe Gamma-LB?  2) ¿Existe Gamma-LBVA?

3) ¿Existe Gamma-LBVF?  4) ¿Existe Gamma-LT?

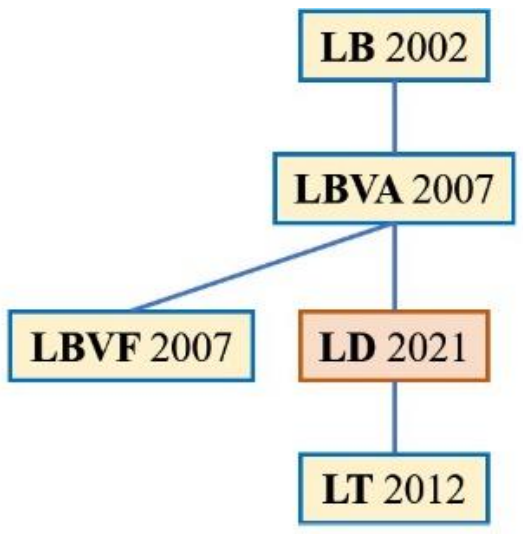

La relación: extensión estricta

# Referencias